\documentclass[11pt]{article}

\title{Structure-preserving integrators for dissipative systems based on reversible-irreversible splitting}

\author{Xiaocheng Shang\footnotemark[1]\ \and Hans Christian \"{O}ttinger\footnotemark[2] }

\date{\today}

\usepackage{cite}

\usepackage{hyperref}
\usepackage{graphicx}

\usepackage{caption}
\DeclareGraphicsExtensions{.eps,.mps,.pdf,.jpg,.png}
\graphicspath{{figures/}{../figures/}}

\usepackage{fancyhdr}

\usepackage{amsmath}
\usepackage{bm}

\usepackage{amsfonts}

\usepackage{dsfont}

\usepackage[titletoc,toc,title]{appendix}

\usepackage[margin=3.0cm]{geometry}

\newcommand{\dd}{{\rm d}}

\newcommand{\bea}{\begin{eqnarray}}
\newcommand{\eea}{\end{eqnarray}}
\newcommand{\be}{\begin{equation}}
\newcommand{\ee}{\end{equation}}

\newcommand{\kB}{k_\mathrm{B}}

\usepackage{mathrsfs}

\usepackage{soul}
\usepackage{color}

\usepackage{algorithm}
\usepackage{algorithmic}

\usepackage{url}
\usepackage{subfigure}

\usepackage{epstopdf}
\DeclareGraphicsExtensions{.eps,.mps,.pdf,.jpg,.png}

\usepackage{latexsym}
\usepackage{mathrsfs}
\usepackage{amssymb}
\usepackage{enumitem}

\usepackage{url}

\usepackage{array}
\newcolumntype{C}[1]{>{\centering\let\newline\\\arraybackslash\hspace{0pt}}m{#1}}

\begin{document}

\maketitle

\renewcommand{\thefootnote}{\fnsymbol{footnote}}

\footnotetext[1]{School of Mathematics, University of Birmingham, Edgbaston, Birmingham, B15 2TT, United Kingdom. Email: x.shang.1@bham.ac.uk }
\footnotetext[2]{Department of Materials, Polymer Physics, ETH Z\"{u}rich, Leopold-Ruzicka-Weg 4, CH-8093 Z\"{u}rich, Switzerland. Email: hco@mat.ethz.ch }

\begin{abstract}
  We study the optimal design of numerical integrators for dissipative systems, for which there exists an underlying thermodynamic structure known as GENERIC (general equation for the nonequilibrium reversible-irreversible coupling). We present a frame-work to construct structure-preserving integrators by splitting the system into reversible and irreversible dynamics. The reversible part, which is often degenerate and reduces to a Hamiltonian form on its symplectic leaves, is solved by using a symplectic method (e.g., Verlet) with degenerate variables being left unchanged, for which an associated modified Hamiltonian (and subsequently a modified energy) in the form of a series expansion can be obtained by using backward error analysis. The modified energy is then used to construct a modified friction matrix associated with the irreversible part in such a way that a modified degeneracy condition is satisfied. The modified irreversible dynamics can be further solved by an explicit midpoint method if not exactly solvable. Our findings are verified by various numerical experiments, demonstrating the superiority of structure-preserving integrators over alternative schemes in terms of not only the accuracy control of both energy conservation and entropy production but also the preservation of the conformal symplectic structure in the case of linearly damped systems.
\end{abstract}

\pagenumbering{arabic}

\section{Introduction}
\label{sec:Introduction}

As an introduction to this article on structure-preserving integrators for dissipative systems, we first summarize the state-of-the-art of the literature and then provide a description of the GENERIC formulation and its properties. The introduction ends with an outline of the article.

\subsection{State-of-the-art in structure-preserving integrators}

In the last few decades, considerable effort has been devoted to developing structure-preserving integrators for Hamiltonian systems. It has been demonstrated that the so-called symplectic integrators, which preserve the symplectic structure, have superior long time behavior compared to their nonsymplectic counterparts, and should be preferred in practice~\cite{Hairer2006,Leimkuhler2005,Leimkuhler2015b}. On the other hand, there has been growing interest in designing appropriate numerical methods for gradient flows~\cite{Stuart1996,Jordan1998,Otto2001,Ambrosio2008,Hairer2014} that respect their underlying properties. In contrast to the symplectic structure, the conformal symplectic structure~\cite{Dressler1988,Moser1994,McLachlan2001,Bhatt2016,Hong2017,Bhatt2017} for Hamiltonian systems that are perturbed by a linear damping (which can be thought of as a special case of the Rayleigh dissipation) has been less studied. It is also worth mentioning that variational integrators~\cite{Kane2000} and specialized Runge--Kutta methods~\cite{Jay2003} have also been used to solve dissipative systems. It turns out that thermodynamically admissible evolution equations for nonequilibrium systems have a more general (including an additional variable known as entropy) and well-defined structure known as GENERIC~\cite{Grmela1997,Oettinger1997,Oettinger2005,Oettinger2018}, which possesses the following distinct features:
\begin{enumerate}
  \item[(i)] conservation of the total energy;
  \item[(ii)] separation of the reversible and irreversible dynamics;
  \item[(iii)] the reversible dynamics preserves a Poisson structure;
  \item[(iv)] entropy production is unaffected by the reversible dynamics;
  \item[(v)] nonnegative entropy production rate.
\end{enumerate}

\subsection{GENERIC formulation}

The GENERIC formulation of the time evolution for nonequilibrium systems is given by
\begin{equation}\label{eq:GENERIC}
  \frac{\dd x}{\dd t} = L \frac{\partial E}{\partial x} + M \frac{\partial S}{\partial x} \, ,
\end{equation}
where $x$ is the set of independent variables required to describe a given nonequilibrium system, $E$ and $S$ represent, respectively, the total energy and entropy as functions of the independent variables $x$, and $L$ and $M$ denote the antisymmetric Poisson matrix and the positive semidefinite (symmetric) friction matrix, respectively. Note that both $L$ and $M$ can also depend on the independent variables $x$ so that the fundamental time evolution equation~\eqref{eq:GENERIC} could be highly nonlinear. We also point out that $\partial/\partial x$ in~\eqref{eq:GENERIC} simply implies the partial derivative although it typically denotes the functional derivative when $x$ is a function/field. Moreover,~\eqref{eq:GENERIC} is supplemented by two degeneracy conditions:
\begin{equation}\label{eq:GENERIC_Degen_1}
  L \frac{\partial S}{\partial x} = 0 \, ,
\end{equation}
and
\begin{equation}\label{eq:GENERIC_Degen_2}
  M \frac{\partial E}{\partial x} = 0 \, .
\end{equation}
Equations~\eqref{eq:GENERIC_Degen_1}--\eqref{eq:GENERIC_Degen_2} indicate the conservation of the entropy by the reversible dynamics (i.e., the $L$ contribution) and the conservation of the total energy in a closed system by the irreversible dynamics (i.e., the $M$ contribution), respectively. Note that ``reversible'' and ``irreversible'' dynamics (in thermodynamics) are simply the names of the two fundamental contributions to the time evolution equation~\eqref{eq:GENERIC}, and should not be confused with similar terms in other subjects. The rank of $M$ has the interpretation of the number of dissipative processes taking place in the system. (See more discussions on the formulation of the GENERIC framework in~\cite{Grmela1997,Oettinger1997,Oettinger2005,Oettinger2018}.)

The usefulness and maturity of the GENERIC framework have been illustrated in a very large number of successful applications in a wide range of areas in Appendix E of~\cite{Oettinger2005} (see also a most recent review of~\cite{Oettinger2017} and references therein). In particular, despite its simple form, we believe that the irreversible dynamics in~\eqref{eq:GENERIC} is the most general form of meaningful irreversible equations in nonequilibrium thermodynamics---it is a belief based on both a very large variety of successful examples and statistical mechanics, so that it can be called knowledge (in particular, as this belief is widely accepted in the nonequilibrium thermodynamics community).

In order to further demonstrate the general properties of $L$ and $M$, the respective Poisson and dissipative brackets are often adopted:
\begin{align}
    \{\mathcal{A},\mathcal{B}\} &= \frac{\partial \mathcal{A}}{\partial x} \cdot L \frac{\partial \mathcal{B}}{\partial x} \, , \label{eq:Poisson_Bracket} \\
    [\mathcal{A},\mathcal{B}] &= \frac{\partial \mathcal{A}}{\partial x} \cdot M \frac{\partial \mathcal{B}}{\partial x} \, , \label{eq:Dissipative_Bracket}
\end{align}
where $\mathcal{A}$ and $\mathcal{B}$ are sufficiently regular (and real-valued) functions of the independent variables $x$. With the help of the two brackets and the chain rule, the time evolution equation of an arbitrary function $\mathcal{A}$ can then be written as
\begin{equation}
  \frac{\dd \mathcal{A}}{\dd t} = \{\mathcal{A},E\} + [\mathcal{A},S] \, .
\end{equation}
More specifically, the Poisson bracket~\eqref{eq:Poisson_Bracket} inherits the antisymmetry of $L$,
\begin{equation}
  \{\mathcal{A},\mathcal{B}\} = - \{\mathcal{B},\mathcal{A}\} \, ,
\end{equation}
and satisfies the Leibniz rule,
\begin{equation}
  \{\mathcal{A}\mathcal{B},\mathcal{C}\} = \mathcal{A} \{\mathcal{B},\mathcal{C}\} + \mathcal{B} \{\mathcal{A},\mathcal{C}\} \, ,
\end{equation}
where $\mathcal{C}$ is another arbitrary sufficiently regular (and real-valued) function of the independent variables $x$. In addition, the Poisson bracket is required to satisfy the Jacobi identity,
\begin{equation}
  \{ \mathcal{A}, \{\mathcal{B},\mathcal{C}\} \} + \{ \mathcal{B}, \{\mathcal{C},\mathcal{A}\} \} + \{ \mathcal{C}, \{\mathcal{A},\mathcal{B}\} \} = 0 \, .
\end{equation}
The dissipative bracket~\eqref{eq:Dissipative_Bracket} inherits the symmetry of $M$,
\begin{equation}
  [\mathcal{A},\mathcal{B}] = [\mathcal{B},\mathcal{A}] \, ,
\end{equation}
and also satisfies the Leibniz rule,
\begin{equation}
  [\mathcal{A}\mathcal{B},\mathcal{C}] = \mathcal{A} [\mathcal{B},\mathcal{C}] + \mathcal{B} [\mathcal{A},\mathcal{C}] \, .
\end{equation}
The positive semidefinite nature of $M$ leads to the nonnegativeness condition
\begin{equation}
  [\mathcal{A},\mathcal{A}] \geq 0 \, ,
\end{equation}
which implies the second law of nonequilibrium thermodynamics (i.e., the entropy production rate is always nonnegative),
\begin{equation}
  \frac{\dd S}{\dd t} = \frac{\partial S}{\partial x} \cdot M \frac{\partial S}{\partial x} = [S,S] \geq 0 \, .
\end{equation}

This article addresses the long-standing challenge of how to preserve the underlying structures when numerically discretizing GENERIC systems in practice. Although in recent years, this topic has attracted increasing attention~\cite{Morrison2017,Portillo2017,Krueger2016,Krueger2011}, to the best of our knowledge, there are no such numerical integrators in the literature. Unlike common approaches that are based on exact energy conservation, we propose in this article a framework to construct structure-preserving integrators for dissipative systems, i.e., GENERIC integrators (also known as metriplectic integrators~\cite{Morrison1986,Morrison1984,Kaufman1984,Grmela1984,Kraus2017} in the mathematical literature), based on splitting the reversible and irreversible dynamics. The topic of structure-preserving integrators for GENERIC/metriplectic systems is the counterpart and generalization of the theory of symplectic integrators for Hamiltonian systems.

\subsection{Outline of the article}

The rest of the article is organized as follows. We give specific definitions of GENERIC integrators and discuss their requirements in numerical discretizations in Section~\ref{sec:Definitions_GENERIC}. In Section~\ref{sec:Split_GENERIC}, we propose a framework to construct split GENERIC integrators based on reversible and irreversible splitting, the generality of the framework is demonstrated in examples of linearly damped systems in Section~\ref{subsec:LDS} as well as in a more challenging (and fully coupled) case of two gas containers exchanging heat and volume in Section~\ref{subsec:TGC}. Section~\ref{sec:Numerical_Experiments} presents various numerical experiments to investigate the performance of the two split \mbox{GENERIC} integrators introduced in this article. Our findings are summarized in Section~\ref{sec:Conclusions}.

\section{Definitions of GENERIC integrators}
\label{sec:Definitions_GENERIC}

In this section, we provide the definitions of \mbox{GENERIC} integrators and discuss their requirements when numerically discretizing a system in practice.

\subsection{Full GENERIC integrators}
\label{subsec:Full_GENERIC}

We recall the definition of (full) GENERIC integrators given in~\cite{Oettinger2018}. Analogous to the definition of symplectic integrators for Hamiltonian dynamics~\cite{Moser1968}, a mapping, $x_{0} \mapsto x_{h}$, is said to be a full GENERIC integrator if it corresponds to a continuous time evolution of a modified GENERIC system
\begin{equation}\label{eq:GENERIC_mod}
  \frac{\dd x}{\dd t} = L \frac{\partial \tilde{E}_{h}}{\partial x} + \tilde{M}_{h} \frac{\partial S}{\partial x} \, ,
\end{equation}
where $\tilde{E}_{h}$ and $\tilde{M}_{h}$ represent the modified energy and friction matrix associated with the integrator, respectively, satisfying a modified degeneracy condition:
\begin{equation}\label{eq:GENERIC_Degen_mod}
  \tilde{M}_{h} \frac{\partial \tilde{E}_{h}}{\partial x} = 0 \, .
\end{equation}
That is, given initial conditions $x(0)=x_{0}$, the analytical solution of~\eqref{eq:GENERIC_mod}, $x(t)$, should agree with what we obtain from the integrator at time $h$, i.e., $x(h)=x_{h}$. A full GENERIC integrator $x \mapsto x_{h}$, which can be thought of as the formal solution of~\eqref{eq:GENERIC_mod}, possesses the following structure:
\begin{equation}
  x_{h} = \exp \left\{ h \left( L \frac{\partial \tilde{E}_{h}}{\partial x} + \tilde{M}_{h} \frac{\partial S}{\partial x} \right) \cdot \frac{\partial}{\partial x} \right\} x \, .
\end{equation}
Similar to symplectic integrators for Hamiltonian dynamics, the modified energy, $\tilde{E}_{h}$, is strictly conserved by a GENERIC integrator. The physical energy $E$ is expected to remain close to the modified energy, $\tilde{E}_{h}$, even for long
integration periods. Additionally, the modified friction matrix, $\tilde{M}_{h}$, should not introduce any additional dissipative processes not present in the original matrix $M$. We point out that full GENERIC integrators may only be available in special cases, for instance, a full GENERIC integrator in the case of a damped harmonic oscillator, where analytical solutions of the GENERIC system can be obtained, was proposed and discussed in~\cite{Oettinger2018}. However, it should be noted that it is highly unlikely that analytical solutions would be available for general GENERIC systems. (Nevertheless, it might be eventually possible to recognise a full GENERIC integrator without exact solutions.) Therefore, in what follows we introduce a framework to construct ``split'' GENERIC integrators.

\subsection{Split GENERIC integrators}
\label{subsec:Split_GENERIC}

Inspired by recent developments on splitting methods~\cite{Leimkuhler2013,Leimkuhler2013a,Leimkuhler2015b,Leimkuhler2013c,Abdulle2014,Leimkuhler2015,Leimkuhler2016a,Leimkuhler2015a,Shang2017}, we consider to split the reversible and irreversible parts of the GENERIC system in such a way that the reversible dynamics, which is often degenerate but possesses a Hamiltonian form on its symplectic leaves, can be integrated by using a symplectic method (e.g., Verlet) with degenerate variables being left unchanged, while the irreversible part (gradient flow) can be solved in such a way that as many structure elements as possible can be preserved (see more references on the challenging task of structure preservation on manifolds in~\cite{Stuart1996,Jordan1998,Otto2001,Ambrosio2008,Hairer2014,Matthes2014,Matthes2019}).

An interesting question for the split GENERIC integrators is: under what conditions do a modified energy and an associated friction matrix, satisfying the modified degeneracy condition~\eqref{eq:GENERIC_Degen_mod}, exist? If they exist, how much do we know about their respective forms? GENERIC integrators share some common features of GENERIC systems discussed at the beginning of this article, which can also be thought of as the requirements for GENERIC integrators. Denoting the Jacobian matrix of the independent variables $x$ as $\Omega$, we have
\begin{enumerate}
  \item[(i)] preservation of the Poisson structure for the reversible dynamics: $\Omega(x_{0}) L(x_{0}) \Omega^{\mathrm{T}}(x_{0}) = L(x_{h})$;
  \item[(ii)] nonnegative entropy production rate: $S(x_{h}) \geq S(x_{0})$;
  \item[(iii)] the modified degeneracy condition~\eqref{eq:GENERIC_Degen_mod} is satisfied with the other~\eqref{eq:GENERIC_Degen_1} being unchanged;
  \item[(iv)] preservation of the rank of the friction matrix: $\mathrm{rank}(\tilde{M}_{h}) = \mathrm{rank}(M)$.
\end{enumerate}
Note that the satisfaction of the modified degeneracy condition~\eqref{eq:GENERIC_Degen_mod} may be based on a truncated modified energy as discussed in Section~\ref{subsubsec:mYBABY}. As pointed out in~\cite{Quispel2008}, it has been proved in~\cite{Zhong1988} that there cannot exist an integrator for ``non-integrable'' Hamiltonian dynamics that preserves both the symplectic (Poisson) structure and the energy (Hamiltonian). In fact, it has been discussed in~\cite{Simo1992} that the preservation of either property has its advantages and disadvantages. While previous attempts to construct structure-preserving integrators for dissipative systems have been relying on the exact conservation of energy (i.e., the energy-conserving discrete gradient methods~\cite{Quispel1996,Cohen2011,McLachlan1999}, see more discussions in Section~\ref{subsubsec:ADG}), there is no obvious reason why integrators that preserve the Poisson structure for the reversible dynamics should be ignored.

\section{Construction of split GENERIC integrators based on reversible-irreversible splitting}
\label{sec:Split_GENERIC}

In this section, we discuss the construction of GENERIC integrators based on splitting the reversible and irreversible parts of the system. In order to satisfy the modified degeneracy condition~\eqref{eq:GENERIC_Degen_mod}, we explore the possibility of adjusting the irreversible part using a modified friction matrix that corresponds to a modified energy associated with the symplectic integrator used for the reversible part.

\subsection{Linearly damped systems}
\label{subsec:LDS}

We first consider a linearly damped system that possesses a natural GENERIC structure~\eqref{eq:GENERIC} with independent variables $x=(q,p,S)$, where $q$ and $p$ represent the position and momentum of the particle, respectively, and $S$ is the entropy of the surrounding thermal bath. While $S$ is an independent variable and thus $\partial S/\partial x = (0,0,1)$, the total energy of the GENERIC system is given by
\begin{equation}\label{eq:Energy}
  E(q,p,S) = H(q,p) + TS = \frac{p^{2}}{2m} + U(q) + TS \, ,
\end{equation}
where $H(q,p)$ represents the Hamiltonian of the particle, $U(q)$ denotes the potential energy, and $TS$ is the energy of the thermal bath. Given the antisymmetric Poisson matrix
\begin{equation}
L=
\left(
  \begin{array}{ccc}
     0 & 1 & 0 \\
    -1 & 0 & 0 \\
     0 & 0 & 0 \\
  \end{array}
\right)
\end{equation}
and the positive semidefinite (symmetric) friction matrix 
\begin{equation}
M=
\left(
  \begin{array}{ccc}
     0 & 0 & 0 \\
     0 & \gamma m T & -\gamma p \\
     0 & -\gamma p & \frac{\gamma p^{2}}{m T} \\
  \end{array}
\right)
=
y y^{\mathrm{T}}
\, , \quad
y = \sqrt{\frac{\gamma}{mT}}
\left(
  \begin{array}{c}
     0  \\
     mT  \\
     -p  \\
  \end{array}
\right) \, ,
\end{equation}
where constant parameters $m$, $\gamma$, and $T$ represent the mass of the particle, the damping rate, and the constant temperature of the thermal bath, respectively, the equations of motion of the GENERIC system can be written as
\begin{align}
    \dot{q} &= \frac{p}{m} \, , \label{eq:LDS_1} \\
    \dot{p} &= F(q) - \gamma p \, , \label{eq:LDS_2} \\
    \dot{S} &= \frac{\gamma p^{2}}{mT} \, , \label{eq:LDS_3}
\end{align}
where $F(q)=-U'(q)$ is the conservative force. Note that in this particular case the symplectic leaves are given by the $(q,p)$ subsystem within the reversible dynamics for constant entropy $S$.

\subsubsection{The YBABY method}
\label{subsubsec:YBABY}

Following the discussions in Section~\ref{subsec:Split_GENERIC}, we suggest to split the GENERIC system~\eqref{eq:LDS_1}--\eqref{eq:LDS_3} into reversible and irreversible parts,
\begin{equation}
  \dd \left[ \begin{array}{c} q \\ p \\ S \end{array} \right] =  \underbrace{\left[ \begin{array}{c} \frac{\partial E}{\partial p} \\ 0 \\ 0 \end{array} \right] \dd t }_\mathrm{A} + \underbrace{\left[ \begin{array}{c} 0 \\ - \frac{\partial E}{\partial q} \\ 0 \end{array} \right] \dd t }_\mathrm{B} + \underbrace{\left[ \begin{array}{c} 0 \\ - \gamma p \\ \frac{\gamma p^{2}}{mT} \end{array} \right] \dd t }_\mathrm{Y} \, .
\end{equation}
Moreover, we can always use a symplectic method (e.g., Verlet) for the reversible dynamics on its symplectic leaves (this is possible in the setting of linearly damped systems~\eqref{eq:LDS_1}--\eqref{eq:LDS_3} where $S$ is an independent variable),
\begin{equation}\label{eq:Splitting_Ham}
  \dd \left[ \begin{array}{c} q \\ p \\ S \end{array} \right] =  \underbrace{\left[ \begin{array}{c} \frac{\partial E}{\partial p} \\ 0 \\ 0 \end{array} \right] \dd t }_\mathrm{A} + \underbrace{\left[ \begin{array}{c} 0 \\ - \frac{\partial E}{\partial q} \\ 0 \end{array} \right] \dd t }_\mathrm{B} \, ,
\end{equation}
while for linearly damped systems with the total energy~\eqref{eq:Energy} the irreversible dynamics
\begin{equation}\label{eq:Splitting_Y}
  \dd \left[ \begin{array}{c} q \\ p \\ S \end{array} \right] =  \underbrace{\left[ \begin{array}{c} 0 \\ - \gamma p \\ \frac{\gamma p^{2}}{mT} \end{array} \right] \dd t }_\mathrm{Y} \, ,
\end{equation}
is exactly solvable
(with $q$ being left unchanged)
\begin{align}
    p_{h} &= \exp\left( -\gamma h \right) p \, , \label{eq:Exact_Sols_P} \\
    S_{h} &= S + \frac{\gamma p^{2}}{mT} \int^{h}_{0} \exp\left( -2 \gamma t \right) \, \dd t = S + \frac{p^{2}}{2mT} \left[ 1 - \exp\left( -2 \gamma h \right) \right] \, . \label{eq:Exact_Sols_S}
\end{align}
Therefore, we can apply the Verlet method to integrate the reversible part~\eqref{eq:Splitting_Ham}
\begin{equation}\label{eq:Operator_Verlet}
  e^{h \hat{\mathcal{L}}_\mathrm{Verlet} } = e^{ \frac{h}{2} \mathcal{L}_\mathrm{B} } e^{ h \mathcal{L}_\mathrm{A} } e^{ \frac{h}{2} \mathcal{L}_\mathrm{B} } \, ,
\end{equation}
and then further split the exact solver~\eqref{eq:Exact_Sols_P}--\eqref{eq:Exact_Sols_S}, $e^{ h \mathcal{L}_\mathrm{Y} }$, for the irreversible part~\eqref{eq:Splitting_Y} to composite a symmetric splitting method, termed ``YBABY'', as
\begin{equation}\label{eq:Operator_YBABY}
  e^{h \hat{\mathcal{L}}_\mathrm{YBABY} } = e^{ \frac{h}{2} \mathcal{L}_\mathrm{Y} } e^{ h \mathcal{L}_\mathrm{Verlet} } e^{ \frac{h}{2} \mathcal{L}_\mathrm{Y} } = e^{ \frac{h}{2} \mathcal{L}_\mathrm{Y} } e^{ \frac{h}{2} \mathcal{L}_\mathrm{B} } e^{ h \mathcal{L}_\mathrm{A} } e^{ \frac{h}{2} \mathcal{L}_\mathrm{B} } e^{ \frac{h}{2} \mathcal{L}_\mathrm{Y} } \, ,
\end{equation}
where $\exp\left(h \mathcal{L}_f \right)$ denotes the phase space propagator associated with the corresponding vector field $f$, with $\mathcal{L}_f$ being the corresponding generator. The generators for each part of the GENERIC system may be written out as follows:
\begin{align}
  \mathcal{L}_\mathrm{A} &= \frac{p}{m} \cdot \nabla_{q} \, , \label{eqn:generator_A} \\
  \mathcal{L}_\mathrm{B} &= F(q) \cdot \nabla_{p} \, , \label{eqn:generator_B} \\
  \mathcal{L}_\mathrm{Y} &= - \gamma p \cdot \nabla_{p} + \frac{\gamma p^{2}}{mT} \cdot \nabla_{S} \, . \label{eqn:generator_Y}
\end{align}
Thus, the generator for the GENERIC system can be written as $\mathcal{L}_\mathrm{GENERIC} = \mathcal{L}_\mathrm{A} + \mathcal{L}_\mathrm{B} + \mathcal{L}_\mathrm{Y}$.
The integration steps of the YBABY method read:
\begin{align}
    p^{n+1/4}   &= \exp\left( - \gamma h/2  \right) p^{n} \, , \label{eq:YBABY_1} \\
    S^{n+1/2}   &= S^{n} + \left[p^{n}\right]^{2} \left[ 1 - \exp\left( - \gamma h \right) \right] / (2mT)  \, , \label{eq:YBABY_2} \\
    p^{n+2/4} &= p^{n+1/4} + (h/2) F(q^{n}) \, , \label{eq:YBABY_3} \\
    q^{n+1}   &= q^{n} + h m^{-1} p^{n+2/4} \, , \label{eq:YBABY_4} \\
    p^{n+3/4} &= p^{n+2/4} + (h/2) F(q^{n+1}) \, , \label{eq:YBABY_5} \\
    p^{n+1}   &= \exp\left( - \gamma h/2 \right) p^{n+3/4} \, , \label{eq:YBABY_6} \\
    S^{n+1}   &= S^{n+1/2} + \left[p^{n+3/4}\right]^{2} \left[ 1 - \exp\left( - \gamma h \right) \right] / (2mT)  \, . \label{eq:YBABY_7}
\end{align}
The order of convergence of a splitting method can be determined by using the Baker--Campbell--Hausdorff formula~\cite{Hairer2006,Leimkuhler2005,Leimkuhler2015b}. For general operators $A$ and $B$, we have
\begin{equation}
  e^{hA}e^{hB} = e^{hZ_{1}} \, ,
\end{equation}
where
\begin{equation}
  Z_{1} = A + B + \frac{h}{2}\langle A, B\rangle + O(h^{2}) \, ,
\end{equation}
with $\langle A, B\rangle=AB-BA$ being the commutator. Subsequently, we can work out
\begin{equation}
  e^{\frac{h}{2}B}e^{hA}e^{\frac{h}{2}B} = e^{hZ_{2}} \, ,
\end{equation}
where
\begin{equation}
  Z_{2} = A + B + O(h^{2}) \, .
\end{equation}
Therefore, a symmetric splitting typically gives second order convergence whereas a nonsymmetric one is generally first order. One can then obtain the associated operator of the YBABY method:
\begin{equation}
  \hat{\mathcal{L}}_\mathrm{YBABY} = \mathcal{L}_\mathrm{A} + \mathcal{L}_\mathrm{B} + \mathcal{L}_\mathrm{Y} + O(h^{2}) \, ,
\end{equation}
which indicates formally second order convergence for the YBABY method~\eqref{eq:Operator_YBABY}. Note that the order of convergence can also be demonstrated by using the Taylor series expansion for the solutions, but the procedure is often tedious. Note also that in principle higher order methods can also be constructed, as in Hamiltonian dynamics~\cite{Yoshida1990}, by suitably composing the operators.

We would also like to point out that while all three subsystems can be solved exactly in linearly damped systems, in cases where the irreversible part is not exactly solvable (see the example of two gas containers exchanging heat and volume in Section~\ref{subsec:TGC}) it is important to solve the irreversible part by using a numerical method that is at least second order so that an overall second order convergence is expected. Alternatively, one could solve the irreversible part by using a numerical method, which could be first order (e.g., the Euler method), and its adjoint method for half a step each, it can be shown that the resulting YBABY$^{\dagger}$ method is self-adjoint (or \emph{symmetric}) and typically has even order (see more discussions in~\cite{Hairer2006,Leimkuhler2005,Leimkuhler2015b}). However, such a method could become implicit, for instance the adjoint method of the Euler method is the implicit backward Euler method.

In the case of $\gamma=0$, the YBABY method reduces to the Verlet method with degenerate variable $S$ being constant, which is a well-known symplectic method that preserves the Poisson structure for the reversible dynamics~\cite{Hairer2006,Leimkuhler2005,Leimkuhler2015b}. Therefore, in order to guarantee the preservation of the Poisson structure for the reversible dynamics, in what follows we will apply the Verlet method for the reversible part, unless otherwise stated.


For linearly damped systems, it has been demonstrated in~\cite{Bhatt2016,Bhatt2017} that numerical methods that preserve the underlying ``conformal symplectic'' structure~\cite{McLachlan2001} are advantageous over alternative schemes. Moreover, high order conformal symplectic and ergodic schemes for stochastic Langevin equation have also been investigated~\cite{Hong2017}.
\par\noindent {\bf Definition 3.1.} \emph{ A numerical method is said to be conformal symplectic if the symplectic two form decays exponentially with a constant decay rate, i.e.,
\begin{equation}\label{eq:Conformal_Symplectic}
  \dd q_{h} \wedge \dd p_{h} = e^{-Kh} \dd q \wedge \dd p \, ,
\end{equation}
where $\wedge$ represents the wedge product and $K>0$ is the constant decay rate. Similarly, a numerical method is said to be symplectic if the symplectic two form is preserved, i.e.,
\begin{equation}
  \dd q_{h} \wedge \dd p_{h} = \dd q \wedge \dd p \, .
\end{equation}
}
We point out that if the prefactor in front of $\dd q \wedge \dd p$ is initially not in an exponential form, we can always rewrite it into an exponential form as long as it is a constant value between zero and one. Following~\cite{Hong2017,Bhatt2016}, we can show that the YBABY method~\eqref{eq:Operator_YBABY} is conformal symplectic:
\begin{equation}\label{eq:Conformal_Symplectic_YBABY}
\begin{aligned}
  \dd q^{n+1} \wedge \dd p^{n+1} &= e^{-\gamma h/2} \dd q^{n+1} \wedge \dd p^{n+2/4} \, , \\
  &= e^{-\gamma h/2} \dd q^{n} \wedge \dd p^{n+2/4} \, , \\
  &= e^{-\gamma h} \dd q^{n} \wedge \dd p^{n} \, .
\end{aligned}
\end{equation}
in which case the decay rate is the physical damping rate, i.e., $K=\gamma$.

We have so far verified the second order convergence for the YBABY method, and its preservation of the Poisson structure for the reversible dynamics as well as the conformal symplecticity.
However, it is unclear under what conditions there exist a modified energy and an associated friction matrix as in~\eqref{eq:GENERIC_Degen_mod}. To this end, in what follows we modify the irreversible part of the system as discussed at the beginning of this section.

\subsubsection{The mYBABY method}
\label{subsubsec:mYBABY}

It is well known that if a symplectic method is used for the reversible dynamics~\eqref{eq:Splitting_Ham}, there exists a modified Hamiltonian, $\tilde{H}_{h}$, in the form of a (typically infinite) series expansion obtained by using backward error analysis~\cite{Reich1999}, which is exactly preserved by the symplectic integrator~\cite{Hairer2006,Leimkuhler2005,Leimkuhler2015b}. In the example of the Verlet method, the modified Hamiltonian is given by
\begin{equation}\label{eq:Hamiltonian_mod}
  \tilde{H}_{h} = \frac{p^{2}}{2m} + U(q) + h^{2} \left( \frac{U''(q)p^{2}}{12m^{2}} - \frac{\left[U'(q)\right]^{2}}{24m} \right) + O(h^{4}) \, .
\end{equation}
In order to identify a modified energy conserved by a GENERIC integrator, we can replace the original energy $E$~\eqref{eq:Energy} by a modified energy, $\tilde{E}_{h}=\tilde{H}_{h}+TS$,
and then try to explore whether we can construct an associated friction matrix, $\tilde{M}_{h}$, in such a way that the modified degeneracy condition~\eqref{eq:GENERIC_Degen_mod} is satisfied. However, it is unlikely that we can find such a friction matrix due to the infiniteness of the series
expansion (and often complicated higher order terms) in the modified energy. Nevertheless, we can truncate the series expansion of the modified energy to certain order in practice, which will introduce some perturbations to the modified energy. For instance, we can use the Verlet method for the reversible part, and then truncate the modified energy up to second order, introducing a perturbation of order four to the modified energy, to obtain
\begin{equation}\label{eq:Energy_mod}
  \tilde{E}_{h} = \frac{p^{2}}{2m} + U(q) + TS + h^{2} \left( \frac{U''(q)p^{2}}{12m^{2}} - \frac{\left[U'(q)\right]^{2}}{24m} \right) \, .
\end{equation}
Subsequently, we can construct the associated modified friction matrix in the fashion of backward error analysis~\cite{Reich1999,Leimkuhler2005,Hairer2006}:
\begin{equation}
  \tilde{M}_{h} = \tilde{y}_{h} \tilde{y}^{\mathrm{T}}_{h} \, ,
\end{equation}
where $\tilde{y}_{h}$ is assumed to be a truncated series expansion up to second order with $y_{i}=[0, a_{i}, b_{i}]^{\mathrm{T}}, i=1,2$:
\begin{equation}
\begin{aligned}
\tilde{y}_{h} = y + h y_{1} + h^{2} y_{2} = \sqrt{ \frac{\gamma}{mT} }
\left(
  \begin{array}{c}
     0 \\
     mT + h a_{1} + h^{2} a_{2} \\
     -p + h b_{1} + h^{2} b_{2} \\
  \end{array}
\right) \, .
\end{aligned}
\end{equation}
In order to satisfy the modified degeneracy condition
\begin{equation}
  \tilde{M}_{h} \frac{\partial \tilde{E}_{h}}{\partial x} = \tilde{y}_{h} \tilde{y}^{\mathrm{T}}_{h} \frac{\partial \tilde{E}_{h}}{\partial x} = 0 \, ,
\end{equation}
which leads to
\begin{equation}
  \tilde{y}^{\mathrm{T}}_{h} \frac{\partial \tilde{E}_{h}}{\partial x} = 0 \, ,
\end{equation}
the following condition has to be satisfied
\begin{equation}
  \left( mT + h a_{1} + h^{2} a_{2} \right) \frac{\partial \tilde{E}_{h}}{\partial p} + \left( -p + h b_{1} + h^{2} b_{2} \right) \frac{\partial \tilde{E}_{h}}{\partial S} = 0 \, ,
\end{equation}
which has a solution
\begin{equation}
  a_{1} = a_{2} = b_{1} = 0 \, , \quad b_{2} = - \frac{U''(q)p}{6m} \, .
\end{equation}
Thus, the modified friction matrix can be written as
\begin{equation}\label{eq:Modified_Friction_Matrix}
\tilde{M}_{h} =
\left(
  \begin{array}{ccc}
     0 & 0 & 0 \\
     0 & \gamma mT & -\gamma p \alpha(q) \\
     0 & -\gamma p \alpha(q) & \frac{\gamma p^{2}\alpha^{2}(q)}{mT} \\
  \end{array}
\right) \, ,
\end{equation}
where the ``modifying factor'' is given by
\begin{equation}\label{eq:Prefactor}
  \alpha(q) = 1 +   \frac{h^{2}U''(q)}{6m} \, .
\end{equation}
Moreover, the modified friction matrix induces a small (second order) perturbation of the physical entropy production
\begin{equation}
  \frac{\dd S}{\dd t} = \frac{\partial S}{\partial x} \cdot \tilde{M}_{h} \frac{\partial S}{\partial x} = \frac{\gamma p^{2} \alpha^{2}(q)}{mT} \geq 0 \, .
\end{equation}
As a result, the irreversible part, incorporating the modified friction matrix~\eqref{eq:Modified_Friction_Matrix}, becomes
\begin{equation}\label{eq:Splitting_Y_m}
  \dd \left[ \begin{array}{c} q \\ p \\ S \end{array} \right] =  \underbrace{\left[ \begin{array}{c} 0 \\ - \gamma p \alpha(q) \\ \frac{\gamma p^{2}\alpha^{2}(q)}{mT} \end{array} \right] \dd t }_\mathrm{Y_{m}} \, ,
\end{equation}
which can be solved exactly (with $q$ being left unchanged)
\begin{align}
    p_{h} &= \exp\left( -\gamma \alpha(q) h \right) p \, , \label{eq:Exact_Sols_P_m} \\
    S_{h} &= S + \frac{\gamma p^{2}\alpha^{2}(q)}{mT} \int^{h}_{0} \exp\left( -2 \gamma \alpha(q) t \right) \, \dd t = S + \frac{p^{2}\alpha(q)}{2mT} \left[ 1 - \exp\left( -2 \gamma \alpha(q) h \right) \right] \, . \label{eq:Exact_Sols_S_m}
\end{align}
In this case, the generator for the modified irreversible dynamics becomes
\begin{equation}
  \mathcal{L}_\mathrm{Y_{m}} = - \gamma p \alpha(q) \cdot \nabla_{p} + \frac{\gamma p^{2}\alpha^{2}(q)}{mT} \cdot \nabla_{S} \, . \label{eqn:generator_Y_m}
\end{equation}
By replacing the Y piece by $\mathrm{Y_{m}}$ in the YBABY method~\eqref{eq:Operator_YBABY}, we can similarly define a symmetric splitting method, termed ``$\mathrm{Y_{m}BABY_{m}}$'' or ``mYBABY'', as
\begin{equation}\label{eq:Operator_mYBABY}
  e^{h \hat{\mathcal{L}}_\mathrm{mYBABY} } = e^{ \frac{h}{2} \mathcal{L}_\mathrm{Y_{m}} } e^{ h \mathcal{L}_\mathrm{Verlet} } e^{ \frac{h}{2} \mathcal{L}_\mathrm{Y_{m}} } = e^{ \frac{h}{2} \mathcal{L}_\mathrm{Y_{m}} } e^{ \frac{h}{2} \mathcal{L}_\mathrm{B} } e^{ h \mathcal{L}_\mathrm{A} } e^{ \frac{h}{2} \mathcal{L}_\mathrm{B} } e^{ \frac{h}{2} \mathcal{L}_\mathrm{Y_{m}} } \, ,
\end{equation}
where the associated operator can be worked out by applying the Baker--Campbell--Hausdorff formula~\cite{Hairer2006,Leimkuhler2005,Leimkuhler2015b} as
\begin{equation}
  \hat{\mathcal{L}}_\mathrm{mYBABY} = \mathcal{L}_\mathrm{A} + \mathcal{L}_\mathrm{B} + \mathcal{L}_\mathrm{Y_{m}} + O(h^{2}) \, ,
\end{equation}
which indicates formally second order convergence for the mYBABY method~\eqref{eq:Operator_mYBABY}. It can be easily shown that all four requirements listed in Section~\ref{subsec:Split_GENERIC} are satisfied for the mYBABY method. The integration steps of the mYBABY method read:
\begin{align}
    p^{n+1/4}   &= \exp\left( - \gamma \alpha\left(q^{n}\right) h/2  \right) p^{n} \, , \label{eq:mYBABY_1} \\
    S^{n+1/2}   &= S^{n} + \left[p^{n}\right]^{2} \alpha\left(q^{n}\right) \left[ 1 - \exp\left( - \gamma \alpha\left(q^{n}\right) h \right) \right] / (2mT)  \, , \label{eq:mYBABY_2} \\
    p^{n+2/4} &= p^{n+1/4} + (h/2) F(q^{n}) \, , \label{eq:mYBABY_3} \\
    q^{n+1}   &= q^{n} + h m^{-1} p^{n+2/4} \, , \label{eq:mYBABY_4} \\
    p^{n+3/4} &= p^{n+2/4} + (h/2) F(q^{n+1}) \, , \label{eq:mYBABY_5} \\
    p^{n+1}   &= \exp\left( - \gamma \alpha\left(q^{n+1}\right) h/2 \right) p^{n+3/4} \, , \label{eq:mYBABY_6} \\
    S^{n+1}   &= S^{n+1/2} + \left[p^{n+3/4}\right]^{2} \alpha\left(q^{n+1}\right) \left[ 1 - \exp\left( - \gamma \alpha\left(q^{n+1}\right) h \right) \right] / (2mT)  \, . \label{eq:mYBABY_7}
\end{align}
Note that in the case of the ``modifying factor''~\eqref{eq:Prefactor} being unity, the mYBABY method~\eqref{eq:Operator_mYBABY} reduces exactly to the YBABY method~\eqref{eq:Operator_YBABY}.

It can be shown that the truncated energy $\tilde{E}_{h}$~\eqref{eq:Energy_mod} is the truncated modified energy, up to second order, for the mYBABY method~\eqref{eq:Operator_mYBABY}, based on the fact that: (i) the Verlet method for the reversible dynamics preserves $\tilde{E}_{h}$~\eqref{eq:Energy_mod} at second order; (ii) the exact solver for the irreversible dynamics preserves $\tilde{E}_{h}$~\eqref{eq:Energy_mod} exactly. In principle, we could truncate the modified energy $\tilde{E}_{h}$ at higher orders (e.g., fourth, sixth,\dots) than that of~\eqref{eq:Energy_mod}, which would lead to higher orders for the overall methods if the irreversible dynamics can be solved exactly. Moreover, it might be more appropriate to refer those GENERIC integrators that incorporate the truncation of the modified energy to ``pseudo-GENERIC integrators'' (in a sense similar to pseudo-symplectic integrators that preserve the symplectic structure only to certain orders~\cite{Aubry1998}).

%
%
%

It can be further shown that the mYBABY method~\eqref{eq:Operator_mYBABY} preserves the conformal symplectic structure if the Hessian of the potential energy is a constant, i.e., $U''(q)=C$. That is, following~\eqref{eq:Conformal_Symplectic_YBABY}, we have \begin{equation}\label{eq:Conformal_Symplectic_mYBABY}
\begin{aligned}
  \dd q^{n+1} \wedge \dd p^{n+1} &= e^{-\gamma_{\mathrm{m}} h/2} \dd q^{n+1} \wedge \dd p^{n+2/4} \, , \\
  &= e^{-\gamma_{\mathrm{m}} h/2} \dd q^{n} \wedge \dd p^{n+2/4} \, , \\
  &= e^{-\gamma_{\mathrm{m}} h} \dd q^{n} \wedge \dd p^{n} \, ,
\end{aligned}
\end{equation}
where
\begin{equation}\label{eq:Modified_Decay_Rate}
  \gamma_{\mathrm{m}} = \gamma \left( 1 + \frac{h^{2}C}{6m} \right) \, .
\end{equation}
which can be thought of as a modified decay rate compared to the damping rate in the YBABY method~\eqref{eq:Conformal_Symplectic_YBABY}.

Note that the preservation of the conformal symplectic structure is, in the literature, often associated with a decay rate of exactly the damping rate as in the YBABY method. Therefore, we may interpret that the YBABY method preserves the conformal symplectic structure in a ``strong'' sense whereas the mYBABY method preserves the conformal symplectic structure in a ``weak'' sense.

\subsection{Two gas containers exchanging heat and volume}
\label{subsec:TGC}

\begin{figure}[tb]
\centering
\includegraphics[scale=0.1]{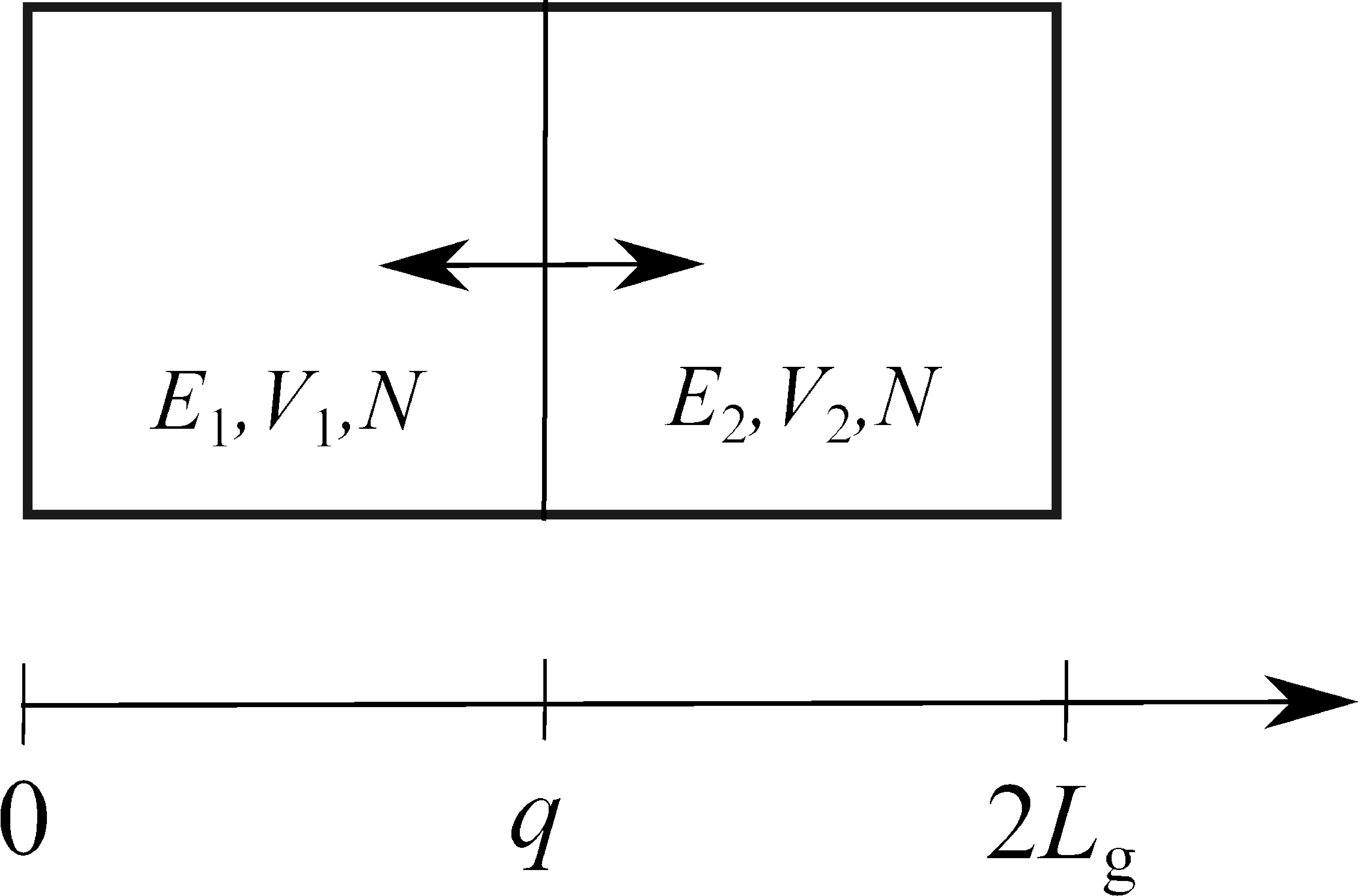}
\caption{\small Schematic illustration of the two gas containers exchanging heat and volume. }
\label{fig:GENERIC_Int_TGC}
\end{figure}

In order to demonstrate the generality of our framework introduced in Section~\ref{subsubsec:mYBABY}, we also consider an example of two (ideal) gas containers exchanging heat and volume (see Fig.~\ref{fig:GENERIC_Int_TGC} and Exercises 3 \& 9 in~\cite{Oettinger2005} for more details) with independent variables $x=(q,p,S_{1},S_{2})$, where $q$ and $p$, respectively, represent the position and momentum of the separating wall of mass $m$, while $S_{1}$ and $S_{2}$ are, respectively, the entropies of the two subsystems. In this case the total energy is given by
\begin{equation}\label{eq:Energy_TGC}
  E(x) = \frac{p^{2}}{2m} + E_{1} + E_{2} \, ,
\end{equation}
where $E_{1}$ and $E_{2}$ are, respectively, the internal energies of the two subsystems with the following relationships to their associated entropies and volumes (i.e., the Sackur--Tetrode equation for ideal gases)
\begin{equation}\label{eq:Entropy_TGC}
  \frac{S_{i}(x)}{N\kB} = \ln \left[ \hat{c} V_{i} \left( E_{i} \right)^{3/2} \right] \, , \quad i=1,2 \, ,
\end{equation}
where $\kB$
is the Boltzmann constant, $\hat{c}$ is another constant that is needed to ensure the argument of the logarithm dimensionless, and it is assumed that the two subsystems contain the same number of particles, $N$. The volumes of the two subsystems are given by
\begin{equation}\label{eq:Volume_TGC}
  V_{1} = q A_{\mathrm{c}} \, , \qquad V_{2} = \left(2L_{\mathrm{g}}-q\right) A_{\mathrm{c}} \, ,
\end{equation}
where $A_{\mathrm{c}}$ is the area of the cross section and $2L_{\mathrm{g}}$ is the length of the container. Given the antisymmetric Poisson matrix
\begin{equation}
L=
\left(
  \begin{array}{cccc}
     0 & 1 & 0 & 0 \\
     -1 & 0 & 0 & 0 \\
     0 & 0 & 0 & 0 \\
     0 & 0 & 0 & 0 \\
  \end{array}
\right)
\end{equation}
and the positive semidefinite (symmetric) friction matrix 
\begin{equation}
M=
\left(
  \begin{array}{cccc}
     0 & 0 & 0 & 0 \\
     0 & 0 & 0 & 0 \\
     0 & 0 & \frac{\alpha}{T^{2}_{1}} & -\frac{\alpha}{T_{1}T_{2}} \\
     0 & 0 & -\frac{\alpha}{T_{1}T_{2}} & \frac{\alpha}{T^{2}_{2}} \\
  \end{array}
\right)
=
y y^{\mathrm{T}}
\, , \quad
y =
\left(
  \begin{array}{c}
     0 \\
     0 \\
     \frac{\sqrt{\alpha}}{T_{1}} \\
   - \frac{\sqrt{\alpha}}{T_{2}} \\
  \end{array}
\right)
\, ,
\end{equation}
where the positive constant parameter $\alpha$ determines the strength of the heat exchange, and $T_{1}$ and $T_{2}$ are, respectively, the temperatures of the two subsystems, related to the associated internal energies by
\begin{equation}\label{eq:Temp_TGC}
  E_{i} = \frac{3}{2} N \kB T_{i} \, , \quad i=1,2 \, ,
\end{equation}
the resulting equations of motion of the GENERIC system can be written as
\begin{align}
    \dot{q} & = \frac{\partial E}{\partial p} = \frac{p}{m} \, , \label{eq:TGC_1} \\
    \dot{p} & = - \frac{\partial E}{\partial q} = \frac{2}{3} \left( \frac{E_{1}}{q} - \frac{E_{2}}{2L_{\mathrm{g}} - q} \right) \, , \label{eq:TGC_2} \\
    \dot{S_{1}} & = \frac{\alpha}{T_{1}} \left( \frac{1}{T_{1}} - \frac{1}{T_{2}} \right) = \frac{9N^{2} \kB^{2} \alpha}{4E_{1}} \left( \frac{1}{E_{1}} - \frac{1}{E_{2}} \right) \, , \label{eq:TGC_3} \\
    \dot{S_{2}} & = - \frac{\alpha}{T_{2}} \left( \frac{1}{T_{1}} - \frac{1}{T_{2}} \right) = - \frac{9N^{2} \kB^{2} \alpha}{4E_{2}} \left( \frac{1}{E_{1}} - \frac{1}{E_{2}} \right) \, . \label{eq:TGC_4}
\end{align}
Although the motion of the wall is assumed to be frictionless (i.e., there is no explicit damping term as in linearly damped systems in Section~\ref{subsec:LDS}), all oscillations of the separating wall have to be damped since they induce (time dependent) temperature differences and thus a heat flux with entropy production. Alternatively, an analysis of the equations~\eqref{eq:TGC_1}--\eqref{eq:TGC_4} linearized around equilibria indicates that the system would always relax to equilibrium. We would also like to point out that, unlike linearly damped systems considered in Section~\ref{subsec:LDS} where the $(q,p)$ dynamics may be viewed as being independent of the entropy, the $(q,p)$ dynamics in this case ``strongly'' depends on the dynamics of $(S_{1},S_{2})$, and vice versa---it is a fully coupled GENERIC system.


As in Section~\ref{subsubsec:YBABY}, we could also split the system~\eqref{eq:TGC_1}--\eqref{eq:TGC_4} into reversible and irreversible parts:
\begin{equation}
  \dd \left[ \begin{array}{c} q \\ p \\ S_{1} \\ S_{2} \end{array} \right] =  \underbrace{\left[ \begin{array}{c} \frac{\partial E}{\partial p} \\ 0 \\ 0 \\ 0 \end{array} \right] \dd t }_\mathrm{A} + \underbrace{\left[ \begin{array}{c} 0 \\  - \frac{\partial E}{\partial q} \\ 0 \\ 0 \end{array} \right] \dd t }_\mathrm{B} + \underbrace{\left[ \begin{array}{c} 0 \\ 0 \\ \frac{\alpha}{T_{1}} \left( \frac{1}{T_{1}} - \frac{1}{T_{2}} \right) \\ - \frac{\alpha}{T_{2}} \left( \frac{1}{T_{1}} - \frac{1}{T_{2}} \right) \end{array} \right] \dd t }_\mathrm{Y} \, ,
\end{equation}
for which we can use the Verlet method for the reversible dynamics, with degenerate variables $S_{1}$ and $S_{2}$ being constants, while a suitable method can be used to solve the irreversible dynamics. We would again like to construct a modified energy and an associated friction matrix as in~\eqref{eq:GENERIC_Degen_mod}. To this end, following the procedures in Section~\ref{subsubsec:mYBABY} we first identify the modified energy associated with the Verlet method used for the reversible dynamics
\begin{equation}\label{eq:Hamiltonian_mod_TGC}
  \tilde{E}_{h} = \frac{p^{2}}{2m} + E_{1} + E_{2} + \frac{h^{2}}{54m}  \left[ \frac{5p^{2}}{m} \left( \frac{E_{1}}{q^{2}} + \frac{E_{2}}{\left(2L_{\mathrm{g}} - q\right)^{2}} \right) - \left( \frac{E_{1}}{q} - \frac{E_{2}}{2L_{\mathrm{g}} - q} \right)^{2} \right] + O(h^{4}) \, ,
\end{equation}
based on which we can subsequently work out the derivatives of the truncated modified energy up to second order:
\begin{equation}
  \frac{\partial \tilde{E}_{h}}{\partial S_{1}} = T_{1} + \frac{h^{2}T_{1}}{54mq} \left[ \frac{5p^{2}}{mq}  - 2 \left( \frac{E_{1}}{q} - \frac{E_{2}}{2L_{\mathrm{g}} - q} \right) \right] \, ,
\end{equation}
and
\begin{equation}
  \frac{\partial \tilde{E}_{h}}{\partial S_{2}} = T_{2} + \frac{h^{2}T_{2}}{54m\left(2L_{\mathrm{g}} - q\right)}  \left[ \frac{5p^{2}}{m\left(2L_{\mathrm{g}} - q\right)} + 2 \left( \frac{E_{1}}{q} - \frac{E_{2}}{2L_{\mathrm{g}} - q} \right) \right] \, .
\end{equation}
We can then construct the associated modified friction matrix in the fashion of backward error analysis as
\begin{equation}\label{eq:Modified_Friction_Matrix_TGC}
\tilde{M}_{h} =
\left(
  \begin{array}{cccc}
     0 & 0 & 0 & 0 \\
     0 & 0 & 0 & 0 \\
     0 & 0 & \frac{\alpha\alpha^{2}_{2}}{T^{2}_{1}} & -\frac{\alpha\alpha_{1}\alpha_{2}}{T_{1}T_{2}} \\
     0 & 0 & -\frac{\alpha\alpha_{1}\alpha_{2}}{T_{1}T_{2}} & \frac{\alpha\alpha^{2}_{1}}{T^{2}_{2}} \\
  \end{array}
\right) = \tilde{y}_{h} \tilde{y}^{\mathrm{T}}_{h}
\, , \quad
\tilde{y}_{h} =
\left(
  \begin{array}{c}
    0 \\
    0 \\
    \frac{\alpha_{2}\sqrt{\alpha}}{T_{1}} \\
  - \frac{\alpha_{1}\sqrt{\alpha}}{T_{2}} \\
  \end{array}
\right) \, ,
\end{equation}
where the modifying factors are given by
\begin{equation}\label{eq:Prefactor_1}
  \alpha_{1} = 1 + h^{2}\beta_{1} \, , \quad \beta_{1} = \frac{1}{54mq} \left[ \frac{5p^{2}}{mq}  - 2 \left( \frac{E_{1}}{q} - \frac{E_{2}}{2L_{\mathrm{g}} - q} \right) \right] \, ,
\end{equation}
and
\begin{equation}\label{eq:Prefactor_2}
  \alpha_{2} = 1 + h^{2}\beta_{2} \, , \quad \beta_{2} = \frac{1}{54m\left(2L_{\mathrm{g}} - q\right)} \left[ \frac{5p^{2}}{m\left(2L_{\mathrm{g}} - q\right)} + 2 \left( \frac{E_{1}}{q} - \frac{E_{2}}{2L_{\mathrm{g}} - q} \right) \right] \, ,
\end{equation}
respectively. The modified irreversible part, incorporating the modified friction matrix~\eqref{eq:Modified_Friction_Matrix_TGC}, is now given by
\begin{equation}\label{eq:Splitting_Y_m_TGC}
  \dd \left[ \begin{array}{c} q \\ p \\ S_{1} \\ S_{2} \end{array} \right] =  \underbrace{\left[ \begin{array}{c} 0 \\ 0 \\ \frac{\alpha\alpha_{2}}{T_{1}} \left( \frac{\alpha_{2}}{T_{1}} - \frac{\alpha_{1}}{T_{2}} \right) \\ - \frac{\alpha\alpha_{1}}{T_{2}} \left( \frac{\alpha_{2}}{T_{1}} - \frac{\alpha_{1}}{T_{2}} \right) \end{array} \right] \dd t }_\mathrm{Y_{m}} \, .
\end{equation}
The modified friction matrix again induces a small (second order) perturbation of the physical entropy production
\begin{equation}
  \frac{\dd S}{\dd t} = \frac{\partial S}{\partial x} \cdot \tilde{M}_{h} \frac{\partial S}{\partial x} = \alpha \left( \frac{\alpha_{2}}{T_{1}} - \frac{\alpha_{1}}{T_{2}} \right)^{2} \geq 0 \, .
\end{equation}
Both YBABY and mYBABY methods are similarly defined in this setting as for linearly damped systems in Section~\ref{subsec:LDS}. However, we are unable to solve the modified irreversible dynamics exactly here, thus a second order explicit midpoint method is suggested to approximate the modified irreversible dynamics~\eqref{eq:Splitting_Y_m_TGC} while the Verlet method is still used for the reversible dynamics, with degenerate variables $S_{1}$ and $S_{2}$ being constants. Overall, the two split GENERIC integrators are both expected to be second order.

We would like to point out that in some cases it might be beneficial to replace the explicit midpoint method by alternative (higher order and/or higher accuracy) methods. Moreover, inspired by the subsampling techniques popular in large-scale Bayesian sampling~\cite{Leimkuhler2015a,Shang2015}, it might be computationally highly advantageous (especially in high dimension) to decompose the positive semidefinite modified friction matrix into nonoverlapping principal submatrices (a principal submatrix can be obtained by selecting a subset of rows and the same subset of columns) that are still positive semidefinite. Having avoided directly solving a high dimensional gradient flow, we could instead solve each of the decomposed and much smaller subsystems with a significantly reduced computational overhead (even with high accuracy). A thorough investigation of this direction is beyond the scope of this article, and will be left for future work.

It is also worth mentioning that when the modified irreversible dynamics has to be approximated by certain numerical methods, the truncated modified energy is expected to be preserved in an ``approximation'' sense. A detailed analysis of the effect of the approximation is also beyond the scope of this article, and will be left for future work.

\section{Numerical experiments}
\label{sec:Numerical_Experiments}

In this section, we conduct various numerical experiments to examine the performance of the two split GENERIC integrators introduced in this article.

\subsection{Simulation details}
\label{subsec:Simulation_Details}

In the case of linearly damped systems, we consider one-dimensional examples of a damped harmonic oscillator (i.e., $U(q)=kq^{2}/2$), for which an analytical solution can be obtained~\cite{Oettinger2018}, as well as a damped nonlinear oscillator (i.e., $U(q)=-k\cos(q)$) where the argument of the cosine function should be dimensionless and this is achieved by fixing the unit of length via the initial position $q_{0}$. The equations of motion of both linearly damped systems can be simplified by dimensional analysis~\cite{Venerus2018}. Without loss of generality, in both cases we choose the basic units (mass, time, temperature, and length, respectively) as $m=k=T=1$ and $q_{0}=2$, where the initial position was particularly chosen to demonstrate the nonlinear effects in the damped nonlinear oscillator. Subsequently, the equations of motion of both linearly damped systems involve only the single dimensionless parameter of $\gamma \geq 0$.
Moreover, in both cases we chose $p_{0} = 0$ as more general values of the initial momentum essentially correspond to a shift of the initial time. Since we are more interested in the deviation from the initial entropy than its absolute value, we set the initial entropy to be zero in both cases.

In the other case of two gas containers (where $A_{\mathrm{c}}=L_{\mathrm{g}}^{2}$) exchanging heat and volume, we chose the basic units of mass and length as $m=L_{\mathrm{g}}=1$, respectively. We further set $N\kB=1$, which fixes a characteristic macroscopic unit of entropy, the counterpart of $T=1$ (i.e, the thermodynamic unit) in the previous examples of linearly damped systems. In order to fix the fourth unit of time, $\alpha=0.5$ was chosen so that: (i) the period of the oscillation is of order one; (ii) there are enough oscillations to collect statistical data (larger values of $\alpha$ lead to faster decay of the amplitude of the oscillation). Furthermore, initial conditions of $(q,p,E_{1},E_{2})=(1,2,2,2)$ were used (i.e., the separating wall is initially in the middle of the container with an initial velocity).

In all three cases, the positions appeared to be oscillating with the amplitudes decaying exponentially. The total simulation time $T_{\mathrm{s}}$ in each case was thus chosen so that $t=T_{\mathrm{s}}$ is the time at which the amplitude of the oscillation was reduced to approximately $1/e$ times its initial value.

Denoting $h$ as the integration stepsize and subsequently $\hat{N}=T_{\mathrm{s}}/h$ the number of integration steps, the root-mean-square error (RMSE) of observable $\phi$ is defined as follows:
\begin{equation}\label{eq:RMSE}
  \mathrm{RMSE} (\phi) = \sqrt{ \frac{1}{\hat{N}} \sum^{\hat{N}}_{i=1} \left( \hat{\phi}_{i} - \phi_{i} \right)^{2} } \, ,
\end{equation}
where $\hat{\phi}_{i}$ and $\phi_{i}$ represent the numerical approximation at time $ih$ and its corresponding exact (reference) value, respectively.


In order to demonstrate the superiority of structure-preserving integrators over alternative schemes, we compare the two split GENERIC integrators introduced in this article with the explicit third order Runge--Kutta (RK3) method used also in~\cite{Bhatt2016} as well as the average discrete gradient (ADG) method~\cite{Harten1983}. The choice of the RK3 method is clearly arbitrary, while other methods are typically second order, it serves as a good example of a higher order method that is not structure-preserving.

\subsubsection{The third order Runge--Kutta method}
\label{subsec:RK3}

Rewriting GENERIC systems in a compact form as $\dot{x}(t) = f(t,x)$ with initial conditions $x(0) = x_{0}$, the RK3 method is given by
\begin{equation}
  x_{n+1} = x_{n} + \frac{h}{6} \left( k_{1} + 4 k_{2} + k_{3} \right) \, ,
\end{equation}
where
\begin{align}
    k_{1} &= f( t_{n}, x_{n} ) \, , \\
    k_{2} &= f( t_{n} + h/2, x_{n} + hk_{1}/2 ) \, , \\
    k_{3} &= f( t_{n} + h, x_{n} - hk_{1} + 2hk_{2} ) \, ,
\end{align}
with $t_{n} = nh, n=0,1,2,\dots$. Note that the RK3 method is neither symplectic nor conformal symplectic. Thus, it does not preserve the Poisson structure for the reversible dynamics. It is also worth mentioning that the two split GENERIC integrators introduced in this article at each step typically require only one force calculation, which often dominates the computational cost per step especially for large-scale simulations, whereas three force calculations are needed for the RK3 method.

\subsubsection{The average discrete gradient method}
\label{subsubsec:ADG}

The so-called discrete gradient methods~\cite{Quispel1996,Cohen2011,McLachlan1999}, which are also known as discrete derivative methods~\cite{Gonzalez1996}, have often been used for the time integration of dissipative systems~\cite{Krueger2011,Romero2009,Hesch2011,Gros2009,Krueger2016,Portillo2017,Suzuki2016,Romero2010,Romero2010a}. For instance, they have recently been suggested to temporally discretize the Landau collision operator in an attempt to preserve its metriplectic/GENERIC structure~\cite{Kraus2017}. However, as stated in~\cite{Cieslinski2010,McLachlan2014a}, discrete gradient methods are generally not symplectic for the symplectic leaves and thus the Poisson structure of the reversible dynamics is not preserved. Therefore, those discrete gradient methods do not belong to either of the GENERIC integrators defined in Section~\ref{sec:Definitions_GENERIC}.

Moreover, discrete gradient methods are typically implicit, in which case iterative methods (e.g., Newton's method) are needed to approximate the solutions at each step. Therefore, discrete gradient methods could be considerably more time-consuming than alternative explicit methods depending on not only the stopping criterion for the iterating procedure~\cite{Cieslinski2010} but also the size of the linear system that needs to be solved at each iteration. However, in the special case of a damped harmonic oscillator (i.e., $U(q)=kq^{2}/2$), we can work out the integration steps without the iterating procedure. More precisely, we rewrite the GENERIC system~\eqref{eq:LDS_1}--\eqref{eq:LDS_3} as
\begin{equation}\label{eq:GENERIC_single_generator_1}
  \frac{\dd x}{\dd t} = \mathcal{S}(x) \nabla E(x) =
\left(
  \begin{array}{ccc}
     0 & 1 & 0 \\
     -1 & - \gamma m & 0 \\
     0 & 0 & \frac{\gamma p^{2}}{m T^{2}} \\
  \end{array}
\right)
\left(
  \begin{array}{c}
     kq \\
     \frac{p}{m} \\
     T \\
  \end{array}
\right)
\, ,
\end{equation}
which is discretized by
\begin{equation}
  \frac{ x_{n+1}-x_{n} }{h} = \bar{\mathcal{S}}(x_{n},x_{n+1}) \bar{\nabla} E(x_{n},x_{n+1}) \, ,
\end{equation}
where the matrix $\bar{\mathcal{S}}(x_{n},x_{n+1})$ approaches $\mathcal{S}(x)$ in the limits of $x_{n+1} \rightarrow x_{n}$ and $h \rightarrow 0$, while the discrete gradient $\bar{\nabla} E(x_{n},x_{n+1})$ satisfies the following conditions:
\begin{align}
  \left( x_{n+1}-x_{n} \right) \cdot \bar{\nabla} E(x_{n},x_{n+1}) &= E(x_{n+1}) - E(x_{n}) \, , \\
  \bar{\nabla} E(x_{n},x_{n}) &= \nabla E(x_{n}) \, .
\end{align}
We consider a midpoint discretization for $\bar{\mathcal{S}}$, i.e.,
\begin{equation}
  \bar{\mathcal{S}}(x_{n},x_{n+1}) = \mathcal{S}(x_{n+1/2}) \, , \quad x_{n+1/2} = \frac{x_{n}+x_{n+1}}{2} \, ,
\end{equation}
and the ADG~\cite{Harten1983} for $\bar{\nabla} E$, i.e.,
\begin{equation}
   \bar{\nabla} E(x_{n},x_{n+1}) = \int^{1}_{0} \nabla E\left( (1-\xi)x_{n} + \xi x_{n+1} \right) \, \dd \xi \, ,
\end{equation}
in which case the ADG method reduces to the implicit midpoint method, which is second order and symplectic (for the symplectic leaf with $\gamma=0$)~\cite{Hairer2006,Hernandez2018}
\begin{align}
    q^{n+1} &= q^{n} + h p^{n+1/2}/m \, , \label{eq:DistGrad_1} \\
    p^{n+1} &= p^{n} - hk q^{n+1/2} - h \gamma p^{n+1/2} \, ,  \label{eq:DistGrad_2} \\
    S^{n+1} &= S^{n} + h \gamma \left[ p^{n+1/2} \right]^{2}/(mT)  \label{eq:DistGrad_3} \, ,
\end{align}
where $q^{n+1/2} = \left( q^{n} + q^{n+1} \right)/2$ and $p^{n+1/2} = \left( p^{n} + p^{n+1} \right)/2$. One might be surprised how, with the irreversible dynamics being switched off (i.e., $\gamma=0$), the energy-conserving ADG method (or the implicit midpoint method) can also be symplectic for the symplectic leaves, which seems to ``contradict'' the findings of~\cite{Zhong1988} (see discussions in Section~\ref{subsec:Split_GENERIC}). However, we point out that in the case of a harmonic oscillator, the corresponding Hamiltonian subsystem is in fact integrable, in such a special case the ADG method preserves not only the energy but also the Poisson structure. Moreover, we can easily solve~\eqref{eq:DistGrad_1}--\eqref{eq:DistGrad_2} to obtain
\begin{align}
    q^{n+1} &= \frac{ \left( 4m + 2mh\gamma -h^{2}k \right)q^{n} + 4hp^{n} }{ 4m + 2mh\gamma + h^{2}k } \, , \\
    p^{n+1} &= \frac{ -4mhkq^{n} + \left( 4m - 2mh\gamma -h^{2}k \right)p^{n} }{ 4m + 2mh\gamma + h^{2}k } \, ,
\end{align}
and subsequently (if $4m + h^{2}k > 2mh\gamma$)
\begin{equation}\label{eq:DistGrad_Wedge}
  \dd q^{n+1} \wedge \dd p^{n+1} = e^{ -\gamma_{\mathrm{ADG}} h} \dd q^{n} \wedge \dd p^{n} \, ,
\end{equation}
where
\begin{equation}\label{eq:ADG_Decay_Rate}
  \gamma_{\mathrm{ADG}} = - \frac{1}{h} \ln \left( \frac{ 4m - 2mh\gamma + h^{2}k }{ 4m + 2mh\gamma + h^{2}k } \right) = \gamma + O(h^2) \, .
\end{equation}
We can see from~\eqref{eq:DistGrad_Wedge}--\eqref{eq:ADG_Decay_Rate} that the ``symplectic two form'' of the ADG method decays exponentially with a constant decay rate, thus the ADG method in this special case preserves the conformal symplectic structure in a ``weak'' sense. Furthermore, when the irreversible dynamics is switched off (i.e., $\gamma=0$),~\eqref{eq:DistGrad_Wedge} reduces to $\dd q^{n+1} \wedge \dd p^{n+1} = \dd q^{n} \wedge \dd p^{n}$, which indicates that the ADG method preserves the Poisson structure for the reversible dynamics. We emphasize here that if the force is nonlinear or alternative discrete gradient approximations (e.g., the midpoint discrete gradient~\cite{Gonzalez1996}, which was used in several methods compared in~\cite{Krueger2011}) are used, the preservation of the conformal symplectic structure and the Poisson structure for the reversible dynamics is expected to be violated while the iterating procedure seems to be unavoidable, which could result in a substantial computational overhead.

The ADG method is unsurprisingly implicit in both cases of the damped nonlinear oscillator (i.e., $U(q)=-k\cos(q)$) and two gas containers exchanging heat and volume. While the former is similar to the damped harmonic oscillator case except replacing $kq$ in~\eqref{eq:GENERIC_single_generator_1} by $k\sin(q)$, whose ADG is still analytically integrable, the latter is more involved. To be more precise, we rewrite the GENERIC system~\eqref{eq:TGC_1}--\eqref{eq:TGC_4} as
\begin{equation}\label{eq:GENERIC_single_generator_2}
  \frac{\dd x}{\dd t} = \mathcal{S}(x) \nabla E(x) =
\left(
  \begin{array}{cccc}
     0 & 1 & 0 & 0 \\
    -1 & 0 & 0 & 0 \\
     0 & 0 & \frac{\alpha}{T^{2}_{1}} & -\frac{\alpha}{T_{1}T_{2}} \\
     0 & 0 & -\frac{\alpha}{T_{1}T_{2}} & \frac{\alpha}{T^{2}_{2}} \\
  \end{array}
\right)
\left(
  \begin{array}{c}
     \nabla E_{q}(x) \\
     \frac{p}{m} \\
     1 \\
     1 \\
  \end{array}
\right)
\, ,
\end{equation}
where
\begin{equation}
  \nabla E_{q}(x) = \frac{2}{3} \left( \frac{E_{2}}{2L_{\mathrm{g}} - q} - \frac{E_{1}}{q} \right) = \frac{2\hat{B}}{3} \left[ \left( 2L_{\mathrm{g}} - q \right)^{-\frac{5}{3}} e^{\frac{2S_{2}}{3N\kB}} - q^{-\frac{5}{3}} e^{\frac{2S_{1}}{3N\kB}} \right] \, ,
\end{equation}
with the constant $\hat{B}$ being defined as
\begin{equation}
  \hat{B} = \left( \hat{c} A_{\mathrm{c}} \right)^{-\frac{2}{3}} \, .
\end{equation}
In this case, the ADG for $\bar{\nabla} E_{q}(x_{n},x_{n+1})$ is no longer analytically integrable, and thus approximated by using the trapezoidal rule
\begin{equation}
\begin{aligned}
   \bar{\nabla} E_{q}(x_{n},x_{n+1}) =& \int^{1}_{0} \nabla E_{q} \left( (1-\xi)x_{n} + \xi x_{n+1} \right) \, \dd \xi \\
   \approx& \frac{\hat{B}}{3} \left[ \left( 2L_{\mathrm{g}} - q_{n} \right)^{-\frac{5}{3}} e^{\frac{2S_{2,n}}{3N\kB}} + \left( 2L_{\mathrm{g}} - q_{n+1} \right)^{-\frac{5}{3}} e^{\frac{2S_{2,n+1}}{3N\kB}} \right] \\
   & - \frac{\hat{B}}{3} \left[ q_{n}^{-\frac{5}{3}} e^{\frac{2S_{1,n}}{3N\kB}} + q_{n+1}^{-\frac{5}{3}} e^{\frac{2S_{1,n+1}}{3N\kB}} \right] \, .
\end{aligned}
\end{equation}
As a result of the approximation, the exact total energy conservation of the ADG method is expected to be violated. (Note that one may rewrite the GENERIC system~\eqref{eq:TGC_1}--\eqref{eq:TGC_4} with independent variables $x=(q,p,E_{1},E_{2})$. However, the same issue of violating the exact total energy conservation for the ADG method is still expected.)

\begin{figure}[tb]
\centering
\includegraphics[scale=0.4]{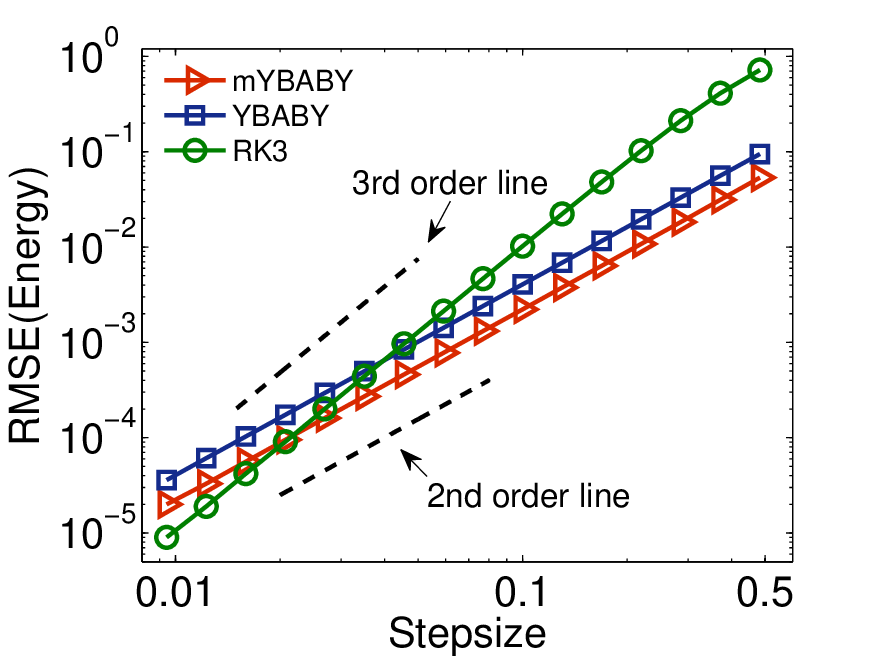}
\includegraphics[scale=0.4]{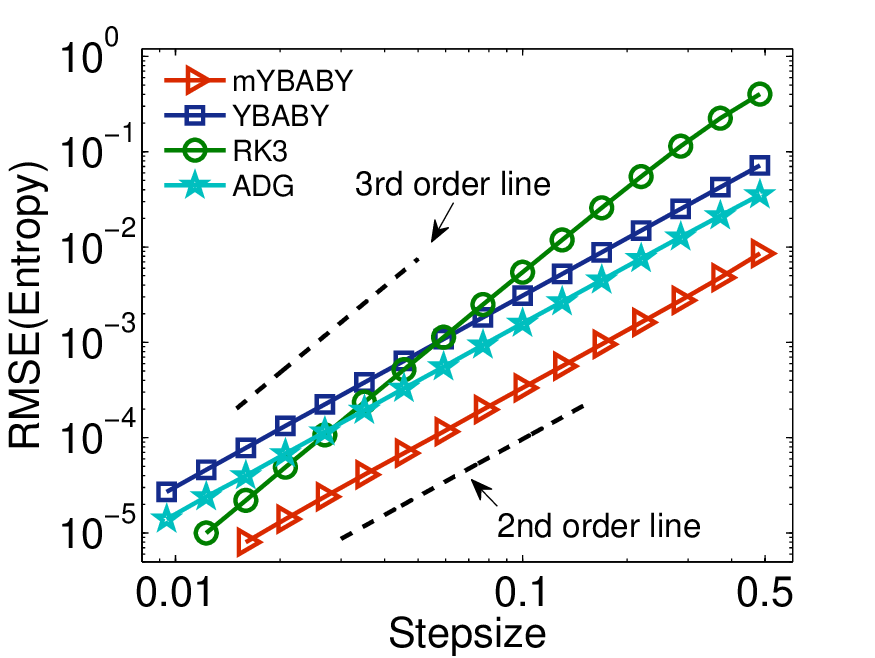}
\caption{\small (Color online) Double logarithmic plot of the root-mean-square error (RMSE)~\eqref{eq:RMSE} in the total energy~\eqref{eq:Energy} (left) and entropy (right) against stepsize by comparing the two split GENERIC integrators introduced in this article with the third order Runge--Kutta (RK3) method and the average discrete gradient (ADG) method, which conserves the total energy exactly (i.e., up to machine precision) and thus is only included in comparisons of the entropy production, with a damping rate of $\gamma=0.01$ and a total simulation time of $T_{\mathrm{s}}=200$ in a standard setting of a damped harmonic oscillator as described in Section~\ref{subsec:Simulation_Details}. The stepsizes tested began at $h=0.0094$ and were increased incrementally by 30\% until around $h=0.5$. Note that, with this set of parameters, the damped harmonic oscillator is ``underdamped'' (i.e., the position of the particle oscillates around zero with the amplitude exponentially decreasing to zero) and the associated period is $T_{\mathrm{p}} \approx 2\pi$. Dashed black lines represent the second and third order convergence as indicated. }
\label{fig:GENERIC_Int_DHO_Energy_Entropy_RMSE}
\end{figure}

\subsection{Damped harmonic oscillator}
\label{subsec:DHO}

We first consider the damped harmonic oscillator (i.e., $U(q)=kq^{2}/2$) example, where the analytical solution is available using the same set of parameters (except $q_{0}$) in~\cite{Oettinger2018}. It is of great importance that numerical approximations of GENERIC systems have (i) a good conservation of the total energy and (ii) a faithful production of the physical entropy, both of which were compared in Fig.~\ref{fig:GENERIC_Int_DHO_Energy_Entropy_RMSE}. We compare the performance of the two split GENERIC integrators with that of the RK3 method and the ADG method. Since the ADG method conserves the total energy exactly (i.e., up to machine precision) in this setting~\cite{Quispel2008}, its results will not be included in comparisons of the energy conservation. According to the dashed order lines, the RK3 method shows third order convergence whereas other methods are all second order as expected. Among the second order methods, the mYBABY method in all the cases we tested outperforms the YBABY method in terms of the RMSE in both quantities. Particularly, in terms of the entropy production, mYBABY is remarkably one order of magnitude more accurate than YBABY, which is slightly outperformed by the ADG method. In both cases, despite its higher computational overhead, the higher order RK3 method is only more accurate than either of the split GENERIC integrators when the stepsize is relatively small, especially for the mYBABY method.

\begin{figure}[tb]
\centering
\includegraphics[scale=0.4]{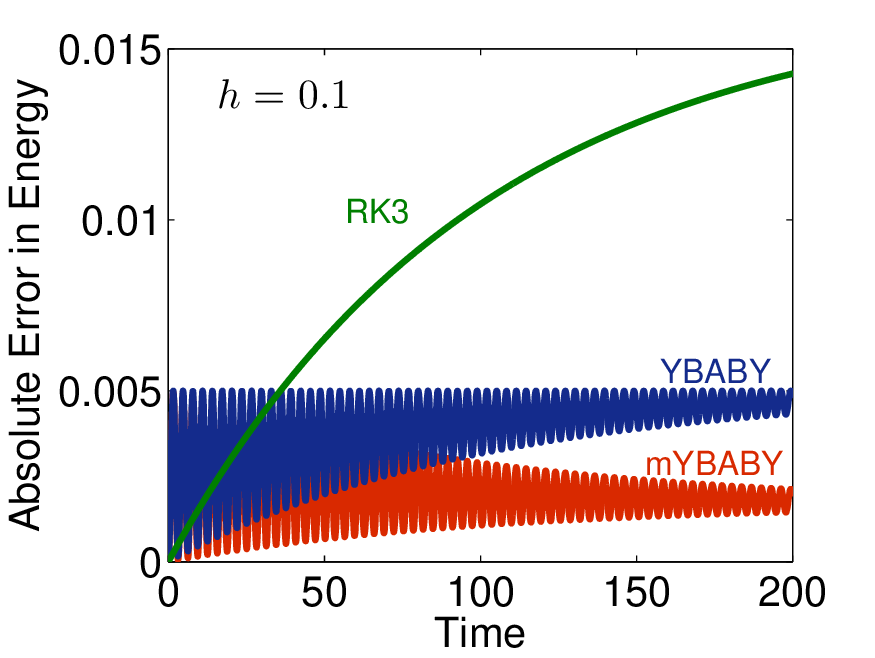}
\includegraphics[scale=0.4]{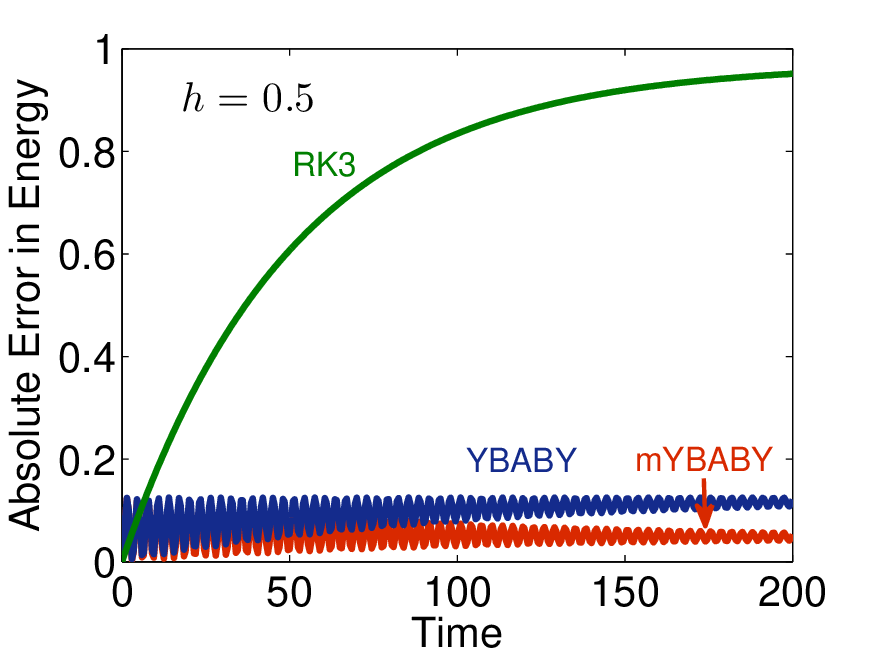}
\caption{\small (Color online) Evolutions of the absolute error in the total energy~\eqref{eq:Energy} obtained from various numerical methods in a standard setting of a damped harmonic oscillator with a damping rate of $\gamma=0.01$, a total simulation time of $T_{\mathrm{s}}=200$, and a fixed stepsize of $h=0.1$ (left) and $h=0.5$ (right).}
\label{fig:GENERIC_Int_DHO_Energy_AE}
\end{figure}

\begin{figure}[tb]
\centering
\includegraphics[scale=0.4]{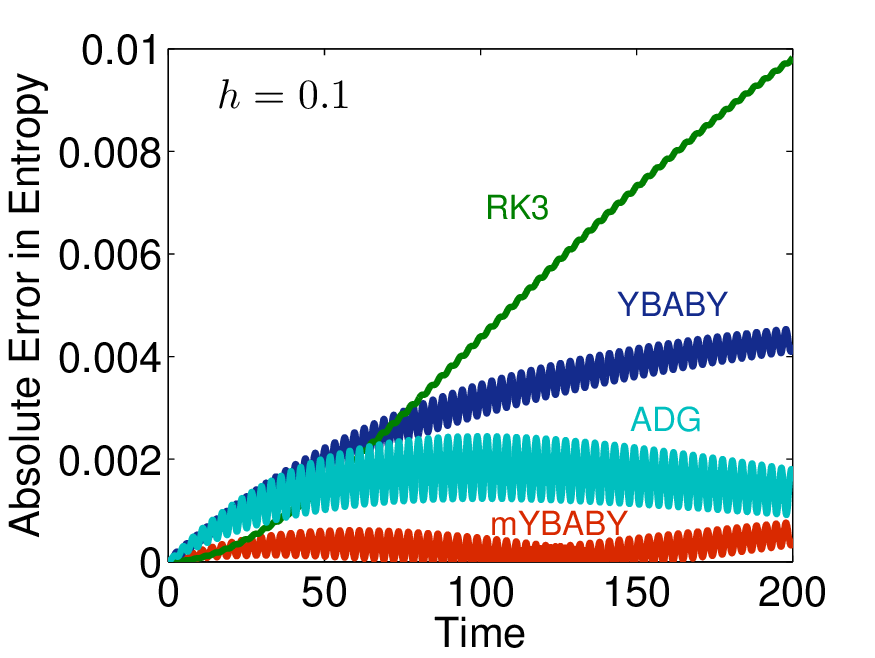}
\includegraphics[scale=0.4]{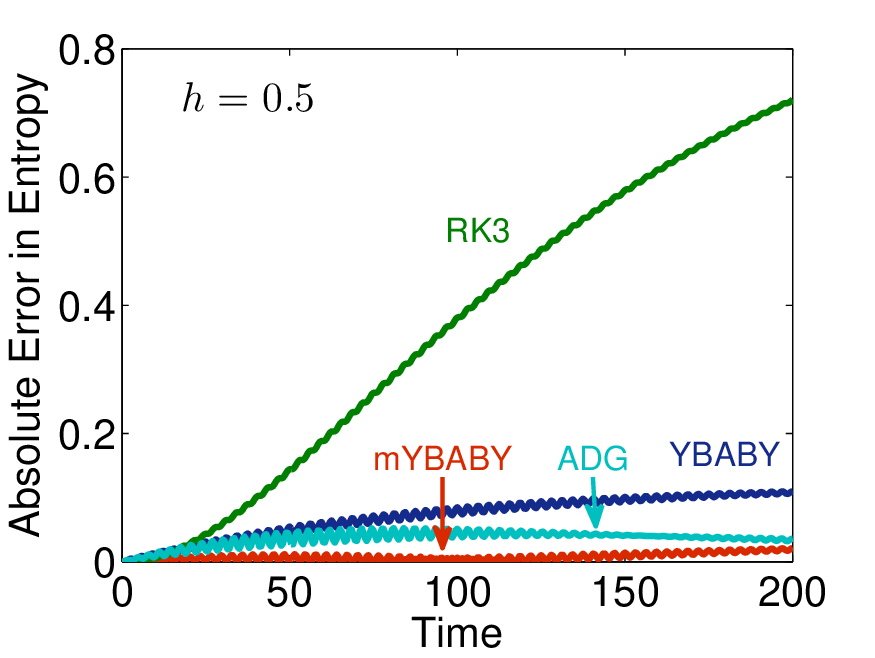}
\caption{\small (Color online) Evolutions of the absolute error in the entropy obtained from various numerical methods in a standard setting of a damped harmonic oscillator with a damping rate of $\gamma=0.01$, a total simulation time of $T_{\mathrm{s}}=200$, and a fixed stepsize of $h=0.1$ (left) and $h=0.5$ (right).}
\label{fig:GENERIC_Int_DHO_Entropy_AE}
\end{figure}

The evolutions of the absolute error in the total energy from various methods were compared (against the exact value of $E_{0}=2$) and plotted in Fig.~\ref{fig:GENERIC_Int_DHO_Energy_AE}. We can see from the figure that, with a stepsize of $h=0.1$ (left panel), the absolute error of the RK3 method rises quickly before eventually settling down, while the absolute errors of the two split GENERIC integrators oscillate strongly with the amplitudes decreasing. Consistent with our findings in Fig.~\ref{fig:GENERIC_Int_DHO_Energy_Entropy_RMSE}, the absolute error of mYBABY is largely smaller than that of YBABY. Nevertheless, both mYBABY and YBABY methods are more accurate than the RK3 method with a relatively large stepsize. The behavior is rather similar with a larger stepsize of $h=0.5$ (right panel) except the magnitude of the error obtained from each method is considerably larger than that with a smaller stepsize.

Figure~\ref{fig:GENERIC_Int_DHO_Entropy_AE} compares the control of the entropy production from a variety of methods. The behavior of the YBABY, mYBABY, and RK3 methods are similar to that of Fig.~\ref{fig:GENERIC_Int_DHO_Energy_AE}. Interestingly, with a stepsize of $h=0.1$ (left panel), the absolute error of the ADG method, while oscillating, initially grows before decreasing while the absolute error of mYBABY, also oscillating, is constantly smaller than that of ADG. The behavior is again very similar with a larger stepsize of $h=0.5$ (right panel) except the magnitude of the errors.

\begin{figure}[tb]
\centering
\includegraphics[scale=0.4]{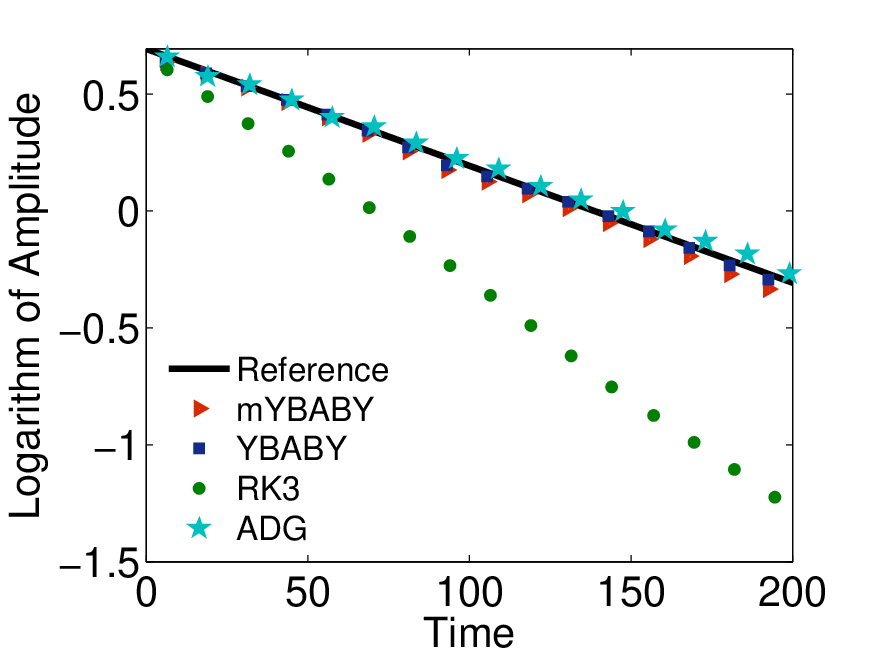}
\caption{\small (Color online) Evolutions of the logarithm of the oscillation amplitude obtained from various numerical methods in a standard setting of a damped harmonic oscillator with a damping rate of $\gamma=0.01$, a total simulation time of $T_{\mathrm{s}}=200$, and a fixed stepsize of $h=0.5$. Note that only the numerical solution of $q$ is presented since it is very similar to that of $p$ with a slight shift. }
\label{fig:GENERIC_Int_DHO_Amplitude_Decay}
\end{figure}

\begin{figure}[tb]
\centering
\includegraphics[scale=0.4]{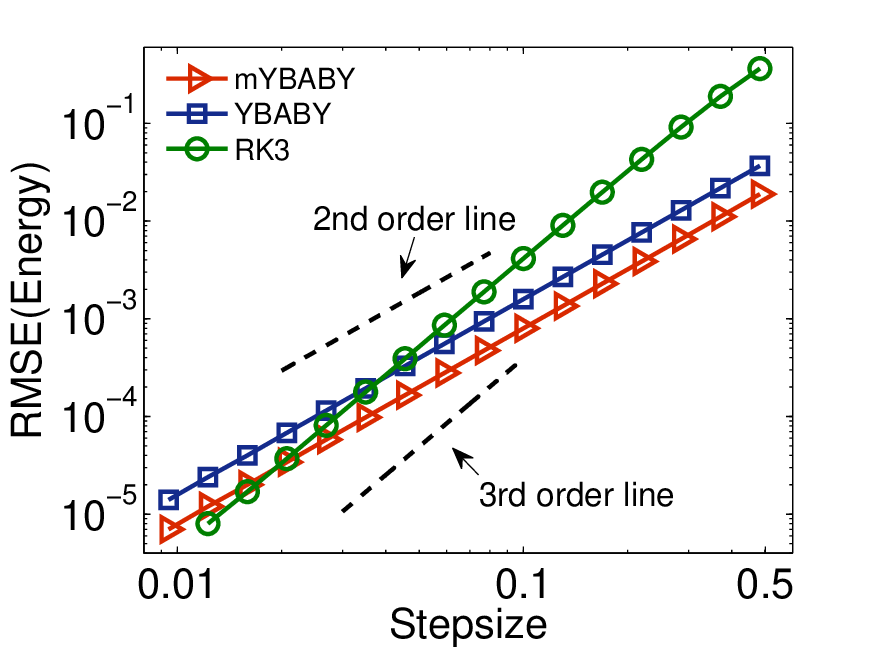}
\includegraphics[scale=0.4]{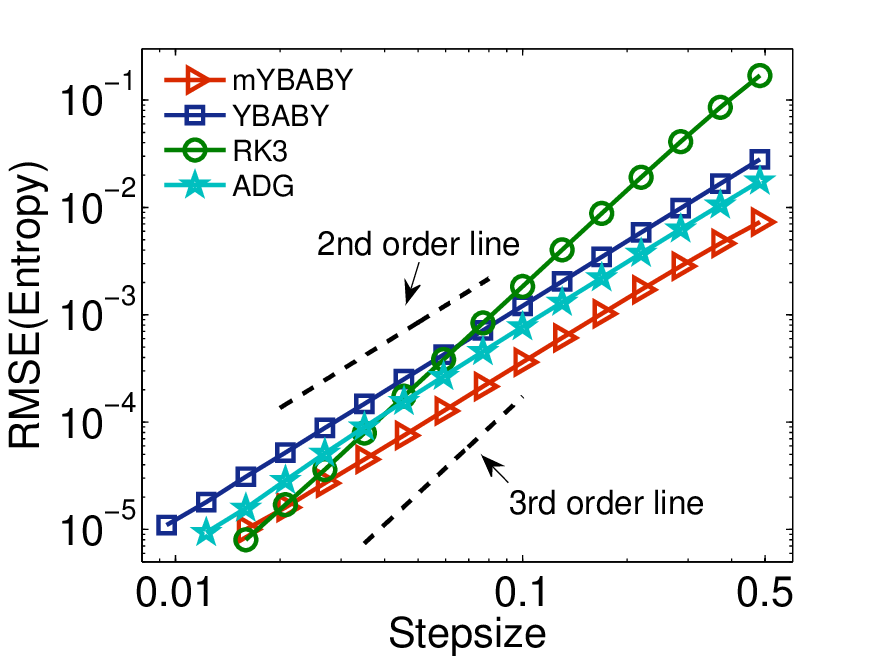}
\caption{\small (Color online) Double logarithmic plot of the root-mean-square error (RMSE)~\eqref{eq:RMSE} in the total energy~\eqref{eq:Energy} (left) and entropy (right) against stepsize by comparing the two split GENERIC integrators introduced in this article with the third order Runge--Kutta (RK3) method and the average discrete gradient (ADG) method with a damping rate of $\gamma=0.01$ and a total simulation time of $T_{\mathrm{s}}=180$ in a standard setting of a damped nonlinear oscillator as described in Section~\ref{subsec:Simulation_Details}. In this case, the position of the damped nonlinear oscillator also oscillates with an associated period of $T_{\mathrm{p}} \approx 8.4$. The format of the plots is the same as in Fig.~\ref{fig:GENERIC_Int_DHO_Energy_Entropy_RMSE}. }
\label{fig:GENERIC_Int_Nonlinear_Energy_Entropy_RMSE}
\end{figure}

\begin{figure}[tb]
\centering
\includegraphics[scale=0.4]{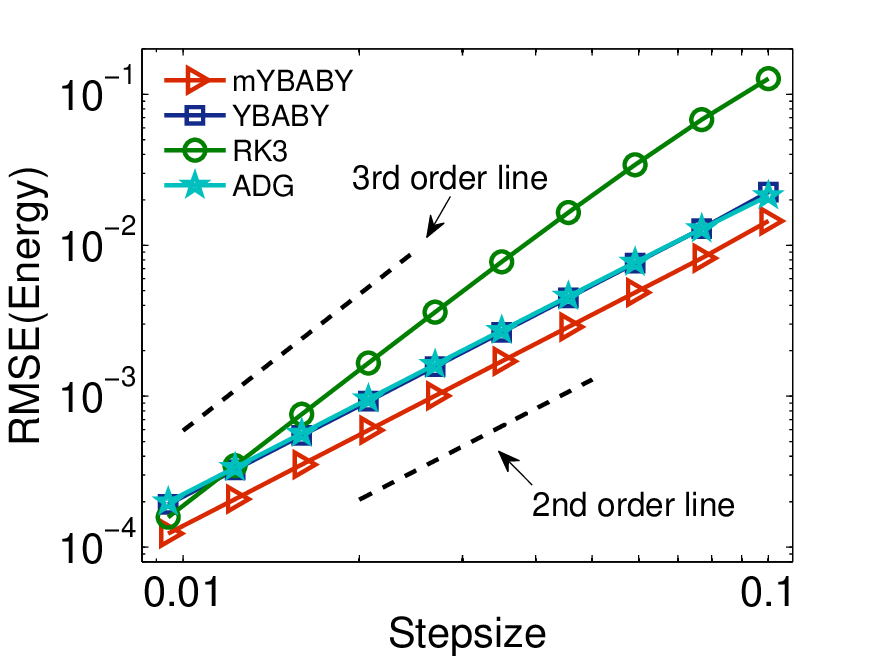}
\includegraphics[scale=0.4]{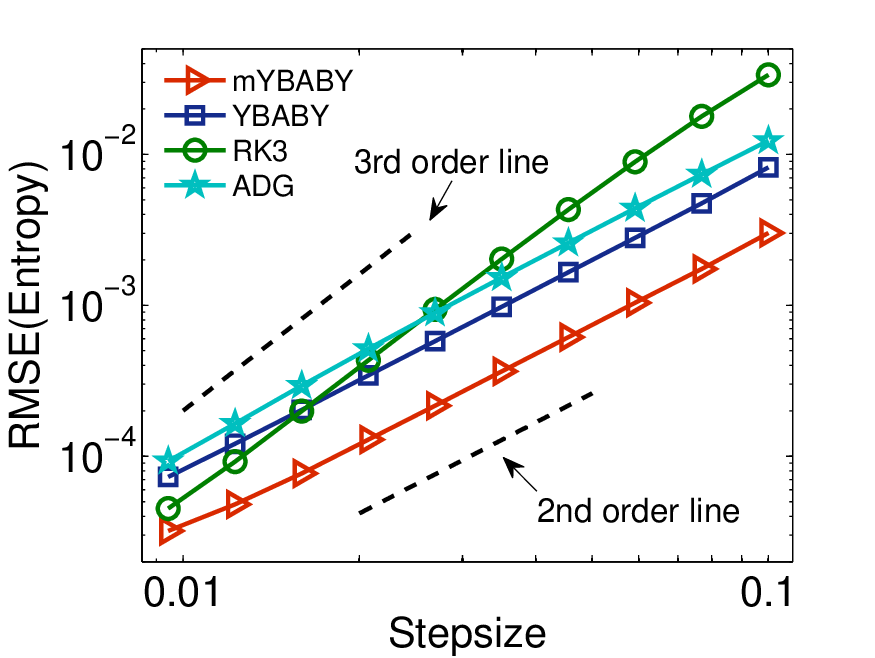}
\caption{\small (Color online) Double logarithmic plot of the root-mean-square error (RMSE)~\eqref{eq:RMSE} in the total energy~\eqref{eq:Energy_TGC} (left) and total entropy (right) against stepsize by comparing the two split GENERIC integrators introduced in this article with the third order Runge--Kutta (RK3) method and the average discrete gradient (ADG) method with a total simulation time of $T_{\mathrm{s}}=30$ in a standard setting of two gas containers exchanging heat and volume as described in Section~\ref{subsec:Simulation_Details}. In this case, the separating wall oscillates around its equilibrium position (i.e., $q=1$) with an associated period of $T_{\mathrm{p}} \approx 2.0$. The format of the plots is the same as in Fig.~\ref{fig:GENERIC_Int_DHO_Energy_Entropy_RMSE}. }
\label{fig:GENERIC_Int_TGC_Energy_Entropy_RMSE}
\end{figure}

We also compare in Fig.~\ref{fig:GENERIC_Int_DHO_Amplitude_Decay} the decay of the oscillation amplitude represented by the ``local maximum'' (in logarithm) of the numerical solution of the position, which characterizes the preservation of the conformal symplectic structure. It can be seen from the figure that while the decay rate of the YBABY method is preserved (almost indistinguishable from the reference decay rate of the damping rate $\gamma$), the RK3 method, which is not conformal symplectic, exhibits a clear drift. It can be also observed (and verified) that both mYBABY and ADG decay at slightly different rates of $\gamma_{\mathrm{m}}$~\eqref{eq:Modified_Decay_Rate} and $\gamma_{\mathrm{ADG}}$~\eqref{eq:ADG_Decay_Rate}, respectively, compared to the reference decay. This indicates that both mYBABY and ADG in this particular case preserve the conformal symplectic structure in a ``weak'' sense (see discussions at the end of Section~\ref{subsubsec:mYBABY}).

\subsection{Damped nonlinear oscillator}
\label{subsec:DNO}

We also investigate the performance of various methods with the damped nonlinear oscillator (i.e., $U(q)=-k\cos(q)$), where the reference solution was obtained by using the RK3 method with a very small stepsize of $h=0.001$. It turns out that the performance of those methods is very similar to that in the case of the damped harmonic oscillator in Section~\ref{subsec:DHO}. Therefore, we only present the results of the accuracy control of both energy conservation and entropy production as in Fig.~\ref{fig:GENERIC_Int_DHO_Energy_Entropy_RMSE}. According to the dashed order lines in Fig.~\ref{fig:GENERIC_Int_Nonlinear_Energy_Entropy_RMSE}, the RK3 method again exhibits third order convergence whereas the other methods appear to be second order. Moreover, we can see from both panels of the figure that the mYBABY method again comfortably outperforms the YBABY method. The RK3 method, with a higher computational overhead, is again in both cases only more accurate than either of the split GENERIC integrators when the stepsize is relatively small. As mentioned in Section~\ref{subsubsec:ADG}, the time-consuming iterating procedure had to be adopted for the ADG method in this nonlinear case. Since the ADG method conserves the total energy up to machine precision, we only include it for comparisons of the entropy production. As can be seen from the right panel of Fig.~\ref{fig:GENERIC_Int_Nonlinear_Energy_Entropy_RMSE} that the ADG method is more accurate than the YBABY method but is outperformed by the mYBABY method despite its higher computational overhead.

\subsection{Two gas containers}

We further examine the performance of various methods in the case of two gas containers exchanging heat and volume described in Section~\ref{subsec:TGC}, where the reference solution was again obtained by using the RK3 method with a very small stepsize of $h=0.001$. While the RK3 method still shows third order convergence, the other methods appear to be second order as expected, according to the dashed order lines in Fig.~\ref{fig:GENERIC_Int_TGC_Energy_Entropy_RMSE}. The performance of the two split GENERIC integrators and the RK3 method is largely similar to that in the previous two examples. More precisely, in both cases the mYBABY method still clearly outperforms the YBABY method, and the two split GENERIC integrators appear to be more accurate than the RK3 method unless the stepsize is relatively small. Unlike the damped nonlinear oscillator example in Section~\ref{subsec:DNO}, the ADG is not analytically integrable and had to be approximated, resulting in errors in the total energy. As a result, we can see from Fig.~\ref{fig:GENERIC_Int_TGC_Energy_Entropy_RMSE} that the performance of the ADG method is almost indistinguishable from that of the YBABY method in terms of the accuracy control of the energy conservation on the left panel, while the former is outperformed by the latter in terms of the accuracy control of the entropy production on the right panel. In both cases, the ADG method is clearly outperformed by the mYBABY method, although the latter only preserves the truncated modified energy in an ``approximation'' sense. Moreover, we would like to point out that while the two split GENERIC integrators are both explicit, the implicit ADG method is computationally much more time-consuming (in this particular case, the evolution of the system was obtained by using the iterative Newton's method at each step in which a linear system associated with a $4\times4$ Jacobian matrix was repeatedly solved).

\section{Conclusions}
\label{sec:Conclusions}

We have given specific definitions of GENERIC integrators that preserve the underlying thermodynamic structures. In order to construct such integrators, we have presented a framework by splitting a GENERIC system into reversible and irreversible parts. The former, which is often degenerate and reduces to a Hamiltonian form on its symplectic leaves, is solved by a symplectic (Verlet) method (with degenerate variables being left unchanged) for which an associated modified Hamiltonian (and subsequently a modified energy) can be obtained by using backward error analysis. The modified energy is subsequently used to construct a modified friction matrix associated with the irreversible part in such a way that the modified degeneracy condition~\eqref{eq:GENERIC_Degen_mod} is satisfied. Following the framework, the mYBABY method has been proposed, which, along with another split GENERIC integrator of the YBABY method, is expected to be second order and typically require only one force calculation at each step. Between the two split GENERIC integrators, we have observed that mYBABY clearly outperforms YBABY in all the cases tested, indicating the importance of satisfying the modified degeneracy condition~\eqref{eq:GENERIC_Degen_mod}.

We have demonstrated by conducting a variety of numerical experiments (including linearly damped systems and two gas containers exchanging heat and volume) that, in terms of the accuracy control of both energy conservation and entropy production, the two split GENERIC integrators (particularly the mYBABY method) are more accurate than the higher order RK3 method unless the stepsizes are relatively small, not to mention the latter requires three force calculations at each step. While the two split GENERIC integrators preserve the conformal symplectic structure for linearly damped systems, RK3 fails and exhibits a clear drift in the decay of the oscillation amplitude of the numerical solutions.

Since the ADG method conserves the total energy up to machine precision, we do not include it in comparisons of the energy conservation for linearly damped systems. It turns out that in both examples of linearly damped systems, the ADG method appears to be more accurate than YBABY, but (despite the use of the time-consuming iterating procedure in the damped nonlinear oscillator case) outperformed by mYBABY in terms of the accuracy control of the entropy production. The ADG is not analytically integrable and had to be approximated in the case of two gas containers exchanging heat and volume, leading to errors in the total energy. As a result of that approximation, the ADG method is as accurate as the YBABY method in terms of the accuracy control of the energy conservation, while the former is outperformed by the latter in terms of the accuracy control of the entropy production. In both cases, although preserving the truncated modified energy in an ``approximation'' sense, the mYBABY method is clearly more accurate than the ADG method. This indicates that in cases where approximations have to be made in the ADG method, it could lose its ``built-in'' advantage of exact conservation of the total energy and be outperformed by alternative methods (especially mYBABY). Moreover, we would like to emphasize again that the implicit ADG method is considerably more time-consuming than the two split GENERIC integrators due to the use of the iterative Newton's method where a linear system associated with a $4\times4$ Jacobian matrix was repeatedly solved at each step. It is anticipated that the computational overhead of discrete gradient methods could be substantially increased for large systems, making them unfavorable compared to explicit structure-preserving integrators (especially mYBABY) in practice.

It is worth mentioning that it might be possible to design GENERIC integrators without an explicit construction of the modified energy $\tilde{E}_{h}$. This is related to the question whether the irreversible dynamics (i.e., a vector field) is ``compatible'' with the canonical transformation associated with the time step $h$. Just as the canonical transformation guarantees that a modified energy $\tilde{E}_{h}$ does exist, there might be a criterion for ``compatibility'' of vector fields with a canonical transformation. The next question would be whether such a ``compatibility'' holds only for the physical entropy or for all possible entropies (which would be the original degeneracy). Alternatively, as symplectic integrators can be obtained most easily from a variational principle~\cite{Bateman1931,Dodin2017,Merker2013,Wendlandt1997}, it might be worth looking at irreversible equations with a variational principle for GENERIC integrators.

\section*{Acknowledgements}

The authors thank Michael Kraus, Martin Kr\"oger, Benedict Leimkuhler, Christian Lubich, Alberto Montefusco, Alexander Ostermann, Jes\'{u}s Mar\'{\i}a Sanz-Serna, Gabriel Stoltz, and Gilles Vilmart for valuable suggestions and comments.

\bibliographystyle{is-abbrv}

\bibliography{refs}

\end{document}